\documentclass[trsc,nonblindrev]{informs3-neutral}

\OneAndAHalfSpacedXI

\usepackage{multirow}

\newlength\shortcolseparator 
\newlength\longcolseparator 

\usepackage[boxed,linesnumbered,noend]{algorithm2e}
\DontPrintSemicolon
\SetAlCapSkip{0.5em}
\SetKwComment{tcp}{\hfill $\triangleright$\ }{}
\SetKwFor{While}{while}{}{end}%
\SetKwFor{For}{for}{}{end}%

\usepackage[hidelinks]{hyperref}
\usepackage[shortlabels]{enumitem} 
\usepackage{microtype}
\usepackage[table]{xcolor}
\usepackage{adjustbox}
\usepackage[caption=false]{subfig}
\usepackage{booktabs, tabularx}
\usepackage{bm} 

\newcommand{\myblue}[1]{#1}

\newcolumntype{M}[1]{>{\centering\arraybackslash\hspace{0pt}}p{#1}}
\newcolumntype{H}{>{\setbox0=\hbox\bgroup}c<{\egroup}@{}}

\usepackage{natbib}
 \bibpunct[, ]{(}{)}{,}{a}{}{,}%
 \def\bibfont{\small}%
 \def\BIBand{and}%

\TheoremsNumberedThrough 
\ECRepeatTheorems

\EquationsNumberedThrough 

\begin{document}


\RUNAUTHOR{Ferraz, Ahmed, Cappart, and Vidal}

\RUNTITLE{Deep Learning for Data-Driven Districting-and-Routing}

\TITLE{Deep Learning for Data-Driven Districting-and-Routing}

\ARTICLEAUTHORS{%
\AUTHOR{Arthur Ferraz}
\AFF{Department of Computer Science, Pontifical Catholic University of Rio de Janeiro, Brazil.}
\AUTHOR{Cheikh Ahmed}
\AFF{
CIRRELT \& SCALE-AI Chair in Data-Driven Supply Chains  \\
Department of Mathematics and Industrial Engineering, Polytechnique Montréal, Canada.
} 
\AUTHOR{Quentin Cappart}
\AFF{
CIRRELT \& Department of Computer and Software Engineering, Polytechnique Montréal, Canada\\
UCLouvain, Louvain-la-Neuve, Belgium.
}
\AUTHOR{Thibaut Vidal}
\AFF{
CIRRELT \& SCALE-AI Chair in Data-Driven Supply Chains  \\
Department of Mathematics and Industrial Engineering, Polytechnique Montréal, Canada.\\
\EMAIL{thibaut.vidal@polymtl.ca}}
} 

\ABSTRACT{%
Districting-and-routing is a strategic problem aiming to aggregate basic geographical units 
(e.g., zip codes) into delivery districts. Its goal is to minimize the expected long-term routing cost of performing deliveries in each district separately. Solving this stochastic problem poses critical challenges since repeatedly evaluating routing costs on a set of scenarios while searching for optimal districts takes considerable time. Consequently, solution approaches usually replace the true cost estimation with continuous cost approximation formulas extending Beardwood-Halton-Hammersley and Daganzo's work. These formulas commit errors that can be magnified during the optimization step. To reconcile speed and solution quality, we introduce a supervised learning and optimization methodology leveraging a graph neural network for delivery-cost estimation. This network is trained to imitate known costs generated on a limited subset of training districts. It is used within an iterated local search procedure to produce high-quality districting plans. Our computational experiments, conducted on five metropolitan areas in the United Kingdom, demonstrate that the graph neural network predicts long-term district cost operations more accurately, and that optimizing over this oracle permits large economic gains (10.12\% on average) over baseline methods that use continuous approximation formulas or shallow neural networks. Finally, we observe that having compact districts alone does not guarantee high-quality solutions and that other learnable geometrical features of the districts play an essential role.
}

\maketitle

\section{Introduction}
\label{sec:introduction}

Districting is the process of partitioning a service region, represented as a collection of basic geographical units, into larger clusters called districts. This practice is ubiquitous in large-scale transportation and last-mile delivery systems for mail delivery \citep{Bruno2021}, home care services \citep{Benzarti2013}, and maintenance services \citep{Garcia-Ayala2016}. A delivery policy in fixed districts has multiple benefits: allowing the separation and the aggregation of the requests in advance before all information is available, reducing the complexity of the task thanks to the decomposition of the routing optimization process, stimulating the familiarity of drivers and thus their efficiency within their respective geographical regions \citep{Zhong2007}, and increasing the satisfaction of customers thanks to a higher familiarity with their drivers \citep{Kovacs2014}.

Districting decisions are strategic and involve significant financial and societal stakes. These decisions typically hold for months or years, whereas operational routes occur on a daily or weekly basis and are subject to variation. Optimizing or even evaluating districting decisions is a very complex task. Demands are uncertain and volatile, and routing cost evaluations typically translate into large-scale vehicle routing problems, which are time-consuming to solve and highly sensitive to the spatial distribution of the requests. Because of these two different classes of decisions and planning horizons, the related \emph{districting-and-routing} problems still represent important challenges \citep{Drexl2015,kalcsics2019}.

An accurate method to estimate costs over a long planning horizon is to consider each district, use a Sample Average Approximation approach (SAA -- \citealt{Verweij2003}) to generate sample demand scenarios (e.g., 10 to 500) within it, and solve a routing problem for each scenario and district. This approach offers the benefit of obtaining a good estimation, but at the cost of extensive computational time due to the large number of districts and scenarios. It is typically applicable to a \emph{fixed} districting solution but becomes impractical within search algorithms for districting problems that require evaluating \emph{numerous} districting solutions (e.g., a local search). 
Therefore, approaches for districting-and-routing regularly rely on continuous approximation formulas to evaluate costs~\citep{Franceschetti2017} without extensive computational overhead. It is, for instance, the case of the Beardwood-Halton-Hammersley (BHH) formula for the traveling salesman problem \citep{Beardwood1959}, which is frequently applied to estimate routing costs through $n$ independently uniformly distributed points in a compact area of size $A$ as $\alpha \sqrt{nA}$, where $\alpha$ is a constant. Though extremely fast, this evaluation approach is far less accurate than a scenario-based approach and may even drive the search toward solutions that have sub-optimal characteristics (see Section~\ref{subsec:impact-on-districting}).

Against this background, we leverage the recent progress in deep learning algorithms, specifically \emph{Graph Neural Networks} \citep[GNN --][]{Scarselli2009,kipf2016semi}, to accurately model expected routing costs and investigate their impact within solution methods for strategic districting-and-routing. Predictions with this approach should be fast enough to support solution-cost evaluations within a local search algorithm, and accurate enough to drive the search toward near-optimal districting solutions.
Geographical areas are divided into geographical units called Basic Units (BU) and represented as a network of BUs with edges linking contiguous areas. Therefore, a solution to the districting problem corresponds to a partition of this graph into connected components representing contiguous districts. Our GNN learns to approximate the routing costs associated with any connected component assimilated to a district. It is trained on a set of districts with known costs obtained using SAA, i.e., by solving TSPs for demand scenarios with the Lin-Kernighan algorithm \citep{LinKernighan1973,Helsgaun2000} and taking their average value as the expected cost.

The trained GNN gives us an oracle that estimates routing costs within districts accurately, as demonstrated by numerical experiments using out-of-sample districts. Furthermore, we must validate whether this performance translates into better strategic districting decisions, since even precise cost estimators may perform badly when integrated within an optimization process due to skew or outliers. To obtain a complete picture, we analyze the impact of using the proposed GNN within search algorithms for districting-and-routing. To that end, we design a prototypical solution approach for districting, which creates initial districts using a mathematical programming method for graph partitioning, and then subsequently applies an Iterated Local Search \citep[ILS --][]{Lourenco2010ILS} to generate better districting plans. The solution approach can use the GNN or any other baseline approach as a cost oracle, so we can conduct extensive computational experiments to compare the impact of different cost estimators on the structure and quality of the final solutions.
Overall, the contributions of this paper are fourfold:
\begin{enumerate}
\item We introduce a learning model based on a GNN architecture for estimating routing costs in a strategic districting problem. This architecture exploits known costs of related districts and the characteristics of the BUs (e.g., population, area, perimeter, density) to provide accurate predictions.
\item We design an ILS that leverages the predicted routing costs to generate high-quality districting plans.
The construction of the initial solution is done by solving an adapted flow formulation for the balanced connected $k$-partition problem, whereas the local search and perturbation procedures are based on BU relocations and exchanges. Our method is generic and can use any cost oracle; it can produce efficient districting plans for practitioners and allows us to gauge the performance of accurate routing-cost estimates in the solution process.
\item We conduct extensive numerical experiments using the data of five metropolitan areas in the United Kingdom (Bristol, Manchester, Leeds, London, and the West Midlands), which contain up to 120 basic geographical units. We compare different cost oracles regarding the quality of routing cost estimates and their performance on the districting task, with districts of different sizes, instances of different scales, and different depot configurations. As seen in these experiments, the proposed GNN-based approach achieves a much smaller prediction error on average compared to other baselines, allowing long-term economic gains of $10.12\%$ on average (and exceeding $20\%$ in some cases) due to better districting decisions.
\item We investigate why cost estimation accuracy has a much larger impact than initially anticipated on the districting-and-routing problem. Indeed, whereas cost estimation accuracy differs by one or two percentage points between methods, their impact on districting-and-routing solution costs is one order of magnitude larger. A finer-grained analysis of the districting solutions shows that classical cost-estimation methods tend to drive the search toward suboptimal solution structures, and that the compactness of the districts is not sufficient to warrant good districting-and-routing solutions in all situations. 
\end{enumerate} 

\section{\myblue{Related Works}}
\label{sec:literature}

\myblue{
Districting problems arise in numerous contexts in which a territory must be partitioned into smaller, connected regions. Comprehensive reviews are available in \citet{Zoltners2005} and \citet{kalcsics2019}, which classify applications ranging from political redistricting to sales, health-care, and logistics planning. In political districting, compactness and fairness are often prioritized \citep[see, e.g.,][]{Horn1993,Webster2013}. In contrast, operational districting, relevant for logistics, delivery, and health-care services, focuses primarily on workload balance, connectivity, and travel cost efficiency. This diversity of objectives has led to a broad spectrum of modeling and algorithmic approaches. To cite a few representative examples related specifically to logistics settings, \cite{Bozkaya2003} introduced a tabu search metaheuristic to solve a multi-criteria districting problem considering compactness and population equality, among other objectives. For home health care districting, \cite{Benzarti2013} introduced and solved different mixed integer linear programming (MILP) formulations, whereas \cite{Garcia-Ayala2016} proposed an arc-based approach for problems with an underlying network structure, including postal delivery, meter readings, winter gritting, road maintenance, and municipal solid waste collection. For large-scale delivery systems, \citet{Bruno2021} analyzed the reorganization of postal collection zones in Bologna to adapt to declining mail volumes, combining spatial and operational constraints. Finally, in the sales and logistics domain, \citet{Lei2015} studied a multi-period districting-and-routing problem with multiple depots and proposed an adaptive large neighborhood search integrating equity, dissimilarity, and profitability objectives.}\\

\noindent
\myblue{
\textbf{Long-term routing cost estimations in districting.}
In distribution contexts, long-term operational efficiency largely depends on the expected delivery costs within each district.  
Accurately evaluating these costs through a Sample Average Approximation (SAA) procedure \citep{Verweij2003}, that is, solving multiple delivery optimization problems such as Traveling Salesman Problems (TSP) or Vehicle Routing Problems (VRP) for sampled demand scenarios, provides reliable estimates but becomes prohibitively time-consuming when embedded in heuristic or metaheuristic frameworks exploring a large number of candidate districts.}

\myblue{
To address this challenge, Continuous Approximation (CA) models \citep{Franceschetti2017,Ansari2018b} have long been employed as analytical surrogates for stochastic evaluations. They are widely used in districting and territory-design studies such as \citet{daganzo1984distance, Novaes2000, Galvao2006, Lei2012} and \citet{Lei2015}. The classical Beardwood-Halton-Hammersley (BHH) theorem \citep{Beardwood1959} expresses the expected tour length through $R$ random points uniformly distributed over an area $A$ as $\alpha \sqrt{RA}$, with $\alpha$ an empirical constant depending on the distance metric.
This relation has inspired numerous refinements. \citet{daganzo1984distance} extended it to vehicle-routing settings by adding a line-haul component proportional to the distance from the depot, while \citet{Chien1992} performed Monte Carlo simulations to determine the best constant in different settings. Some subsequent studies sought to relax the uniformity assumptions: \citet{Figliozzi2007} considered non-uniform demand densities and real urban network effects, and \citet{Cavdar2015b} proposed a distribution-free regression model estimating TSP tour lengths using geometric dispersion and centrality measures. Finally, \citet{Franceschetti2017} extended CA formulations to heterogeneous, time-dependent, and network-based logistics systems.
Although CA methods are extremely fast, they are fundamentally grounded on large-sample approximations and tend to lose accuracy on small or irregular instances. As will be shown later in this paper, their use within districting-and-routing optimization can also steer search heuristics toward geometrical characteristics that are not necessarily aligned with operational efficiency.}\\

\noindent
\myblue{
\textbf{Data-driven estimators of routing costs.} Data-driven research seeks to approximate routing distances from broader data and instance features. These approaches extend the principles of CA formulas, which can be interpreted as linear regression models built on a few spatial descriptors. \citet{Kwon1995} presented an early use of shallow neural networks (with a single hidden layer of three neurons) to estimate the optimal TSP length for uniformly distributed customers in rectangular regions of up to 80 points. More recently, \citet{Akkerman2022} developed linear models with 18 topological features to forecast routing costs in rolling-horizon logistics decisions, while \citet{Varol2024} proposed a three-layer neural network using engineered geometric descriptors of TSP instances and partial solutions generated by heuristics. \citet{Kou2022a} and \citet{Kou2024} introduced computational-statistics estimators based on the distribution of tour lengths obtained from one thousand randomized solutions. However, it is important to note that such estimators offer a practical advantage only when they are faster than solving the underlying problems directly. For medium-sized TSPs (around 100 customers), we observed that an efficient heuristic like Lin-Kernighan (LKH) algorithm \citep{Helsgaun2000} with a short termination criterion tends to be currently faster than the feature-extraction and inference procedures required by some of the aforementioned models.}

\myblue{
Finally, progress has also been made on machine-learning approaches for combinatorial optimization that explicitly \emph{solve} routing problems rather than estimate their cost. Reinforcement-learning and attention-based models have achieved good-quality solutions for medium-sized TSPs \citep[see, e.g.,][]{kool2018attention,joshi2019efficient,kwon2020pomo,sun2023difusco,drakulic2025goal}.
In contrast, other analyses highlighted the limited generalization of these approaches to larger graphs \citep{joshi2022learning}.  
Graph Neural Networks (GNNs) play a central role in this research stream: various architectures \citep[see, e.g.,][]{dai2016discriminative,li2015gated,hamilton2017inductive} encode graph-structured data into latent embeddings that can be leveraged for decision making, and their use in combinatorial optimization has been comprehensively reviewed by \citet{cappart2021combinatorial}.}\\

\myblue{
Overall, despite major advances in deep learning, no advanced machine-learning model has been designed to approximate stochastic district delivery costs from BU-level characteristics, nor has any study examined how estimation strategies affect strategic districting-and-routing outcomes. Existing approaches span from fast but idealized CA formulas to learning-based estimators of deterministic routing costs and, more recently, neural architectures that directly solve routing problems
(both of the latter needing repetitions on different demand scenarios to achieve reliable evaluation accuracy). Yet these methods are either not accurate or fast enough for large-scale districting optimization. This study addresses this gap by introducing a GNN-based methodology to estimate the expected routing cost of each district.
Integrated within an ILS framework, this estimator permits rapid and accurate evaluation of numerous districting configurations, yielding high-quality districting-and-routing solutions.
Furthermore, we provide a systematic analysis of how alternative cost-estimation models influence the structure and efficiency of strategic districting-and-routing decisions.}

\section{Problem Statement}
\label{sec:problem-statement}

We focus on a strategic districting problem encountered when organizing deliveries in a region from a central depot. Consider a geographical region consisting of $n$ geographical units called \textit{Basic Units} (BU). We know the geographical boundaries and the population $\xi_i$ of each BU $i$. A district $d$ is a set of BUs, i.e., $d \, \subseteq \, 2^\Omega$ with $\Omega = \{1,\dots,n\}$. The operational cost of a district is a function $\Phi: 2^{\Omega}\to \mathbb{R}$ representing the expected (i.e., long-term) daily cost of delivering customers in the BUs of this district. Consequently, each district's delivery operations are completely independent.
The districting-and-routing problem studied in this paper then consists of partitioning the region into exactly $k$ districts in such a way that (i) each BU belongs to exactly one district, (ii) the number of BUs inside each district belongs to an admissible range $[n_\textsc{l},n_\textsc{u}]$, (iii) the districts are connected in space, and (iv) the sum of the long-term operational costs of the districts is minimized.\\

\noindent
\textbf{District costs.}
The long-term operational cost $\Phi(d)$ of any district $d$ is calculated as follows. Within each BU~$i$, we assume that a finite set~$X_i(\omega)$ of random demand locations is generated according to the probability distribution $\mathbb{P}_i$ of a spatial Poisson process \citep[see, e.g.,][]{Baddeley2006} with intensity proportional to the population density. Consequently, demand requests are spread uniformly over the geographic area covered by each BU, and BUs with a larger population density have more requests on average per square kilometer.
The cost of any district~$d$ corresponds to the expected distance of the best tour leaving the depot, visiting the customer's locations of all BUs $i \in d$, and returning to the depot. Therefore, the operational cost of~$d$ is defined as:
\begin{equation}
\Phi(d) = \mathbb{E} \left(\min_{\pi} \left\{ \textrm{Dist}(\pi) | \pi \in \textrm{Perm}\left( \bigcup_{i\in d} X_i(\omega)\right)\right\}\right),
\end{equation}
where $\textrm{Perm}(X)$ is the set of all permutations of a set $X$, i.e., each permutation corresponds to a distinct sequence of the elements in $X$, and $\textrm{Dist}(\pi)$ is the travel distance starting from the depot, visiting all locations according to the sequence $\pi$, and returning to the depot.
With this definition, objective function~$\Phi$ is monotonic, i.e., $\Phi(d') \leq \Phi(d)$ if $d' \subseteq d$. However, an exact evaluation of $\Phi$ is intractable as it requires calculating the expected optimal cost of a combinatorial optimization problem with stochastic parameters. To alleviate this issue, we can sample $S$ demand scenarios~$X_{is}$ for $s \in \{1,\dots,S\}$ from~$\mathbb{P}_i$ for each BU $i \in d$, and approximate the cost of the district~as:
\begin{equation}
\Phi_\textsc{saa}(d) = \frac{1}{S} \sum_{s=1}^S \min_{\pi} \left\{ \textrm{Dist}(\pi) | \pi \in \textrm{Perm}\left( \bigcup_{i\in d} X_{is} \right)\right\}.\label{extensive-form}
\end{equation}
Notably, this SAA estimate reduces to solving $|S|$ deterministic TSPs.
\\

\noindent
\textbf{Computational complexity.} Due to the cardinality constraints on the number of BUs per district (between $n_\textsc{l}$ and $n_\textsc{u}$) and the requirement that the districts must remain connected, even identifying an initial feasible solution of the districting-and-routing problem is NP-hard by reduction from the balanced connected $k$-partition problem \citep{Dyer1985}. Moreover, computing the cost of the districts, even considering a single scenario, is another NP-hard subproblem.
\\

\noindent
\textbf{Graph partitioning formulation.} The districting-and-routing problem can be formally cast as a two-stage stochastic graph partitioning problem. Let $G(V, E)$ be an undirected graph, where each vertex $i \in V$ is a BU and edges $e \in E$ represent the contiguity between adjacent BUs (i.e., sharing a border). Let $\mathcal{D} \subset 2^V$ be the set of all possible feasible districts respecting size and connectivity constraints. Then, the problem is formulated as the following integer program:
\begin{align}
 \label{set-part-equation:objective} \text{min} & \sum_{d \in \mathcal{D}} \Phi(d) \lambda_d & \\ 
 \label{set-part-equation:one-bu-per-district} \text{s.t. } & \sum_{d \in \mathcal{D}} b_{id} \lambda_d = 1 & & & i \in V \\
 \label{set-part-equation:k-as-districts} & \sum_{d \in \mathcal{D}} \lambda_d = k & & & \\
 \label{set-part-equation:d-as-binary} & \lambda_d \in \{0,1\} & & & d \in \mathcal{D}.
\end{align}

For each possible district $d$, binary variable $\lambda_d$ takes the value $1$ if and only if $d$ is selected in the solution. Parameter $b_{id} = 1$ if BU $i$ appears in district $d$, and $0$ otherwise. Objective~\eqref{set-part-equation:objective} calculates the total cost of the selected districts. Constraints~\eqref{set-part-equation:one-bu-per-district} ensure that each BU appears in exactly one district, and Constraint~\eqref{set-part-equation:k-as-districts} ensures $k$ districts are formed. This problem can be seen as a two-stage stochastic program in which districts must be formed in the first stage, and one route needs to be selected in each district to visit the revealed customers in the second stage.

Formulation~(\ref{set-part-equation:objective}--\ref{set-part-equation:d-as-binary}) is essentially of descriptive use, as there are two main barriers to its direct solution. First, it includes an exponential number of variables $d \in \mathcal{D}$. Moreover, the calculation of $\Phi(d)$ is very challenging, even when using an approximation as in Equation~\eqref{extensive-form}. For these reasons, this formulation is only practical for obtaining baseline solutions for small problems. For larger cases, more efficient methods are needed. Designing scalable solution methods requires (i) efficient algorithms to estimate district costs (discussed in Section~\ref{sec:neural-network}), and (ii) efficient search strategies for the graph-partitioning problem (discussed in Section~\ref{sec:districting-algorithm}).

\section{Graph Neural Network for District-Cost Estimations}
\label{sec:neural-network}

We focus on the estimation of delivery costs for the districts.
A direct calculation of Equation~\eqref{extensive-form} by solving TSPs on different demand scenarios is fast enough to evaluate the district costs of any given solution of the districting-and-routing problem. However, it is impractically slow within a search method for the districting problem (e.g., a local search) due to the large number of candidate districts considered through the search for move evaluations. To provide fast and accurate estimates for such applications, we explore the option to \emph{learn} an approximation of the delivery costs by supervised learning, more specifically by relying on the considerable methodological progress recently made on deep learning and graph neural networks (GNNs). This approach is described in the remainder of the section, starting from the feature information used for training, the architecture of the network, and the training process.

\subsection{Features}
\label{sec:features}

Our GNN leverages the same undirected graph $G(V, E)$ as the one depicted in Section~\ref{sec:problem-statement}, with the vertices $V$ corresponding to BUs and edges in $E$ corresponding to adjacency relations. For learning the GNN weights, we rely on a training set $\mathcal{D}^\textsc{Train}$ of districts. Each district $d \in \mathcal{D}^\textsc{Train}$ is characterized by a vector of features on each node of the graph and an estimated cost value using SAA as in Equation~\eqref{extensive-form} with $500$ previous demand realization scenarios for each BU.
This estimated cost, also referred to as the \textit{label}, is the value we want the GNN to predict.
We use eight features $\bm{f}_{vd} = (\xi_v, \sqrt{\xi_v}, q_v, a_v, \sqrt{a_v}, \rho_v, \delta_v, e_{vd})^\top \in \mathbb{R}^{8}$ for each vertex~$v$ (associated to a BU) and each district $d$ in the training set. The features used are:
\begin{enumerate}[nosep]
\item the \textbf{population} $\xi_v$ of the BU;
\item the \textbf{sqrt-population} $\sqrt{\xi_v}$ of the BU;
\item the \textbf{area} $a_v$ of the BU;
\item the \textbf{sqrt-area} $\sqrt{a_v}$ of the BU;
\item the \textbf{density} $ \rho_v = \xi_v/a_v$ of the BU;
\item the \textbf{perimeter} $q_v$ of the BU;
\item the \textbf{distance to the depot} $\delta_v$, corresponding to the minimum distance between the depot and any point in the BU;
\item a \textbf{membership variable} $e_{vd}$, taking value $1$ if BU $v$ belongs to $d$ and $0$ otherwise.
\end{enumerate}
The first seven (continuous) features were normalized to $[0,1]$ by dividing by their greatest value.
It is noteworthy that they remain fixed when considering examples generated in the same metropolitan area. Only the last feature changes, according to the current district considered in the sample.
\myblue{
The choice of these features follows well-established CA formulas, notably BHHD, FIG, and SNN (see Section~\ref{sec:literature}). Population and its square root, area and its square root, and the average distance to the depot are standard descriptors in these CA models, to which we additionally include the BU perimeter. Together, these variables summarize the main topological and geographical characteristics of the BUs.} 

\subsection{Architecture of the graph neural network}
\label{sec:GNN-archi}

\begin{figure}[htbp]%
 \centering
 \includegraphics[width=0.99\textwidth]{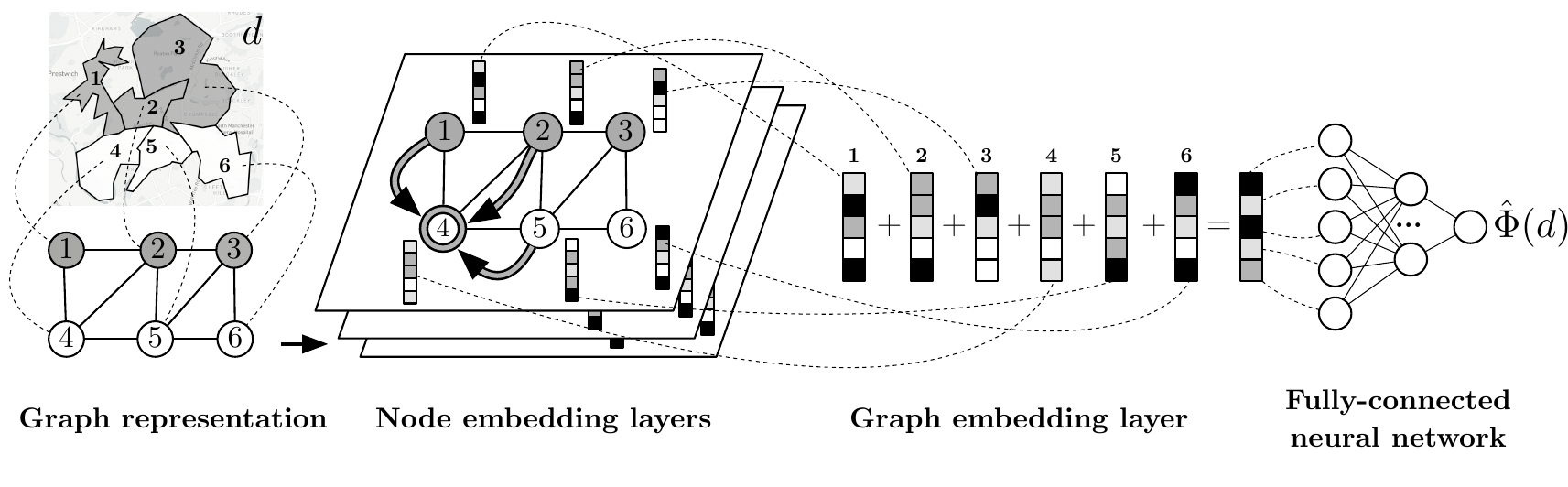}%
 \caption{Neural architecture dedicated to approximate delivery costs.}%
 \label{fig:gnn}%
\end{figure}

The neural network used for the prediction includes three main parts:
\begin{enumerate}
 \item a \textbf{node embedding layer}, taking as input a district sample in the format described previously and whose final output is a latent vector of features for each node of the graph. These outputs are also referred to as \textit{node embedding}. 
 Intuitively, it is a function that converts feature information defined over a graph into a vector representation of the features by aggregating information from neighboring nodes. This aggregation is done several times and corresponds to a layer of the GNN.
 \item a \textbf{graph embedding layer} aggregates each node embedding into a single vector through a non-linear transformation. Intuitively, this vector is a latent representation of the input graph. It is also referred to as the \textit{graph embedding}.
 \item a \textbf{fully-connected neural network} whose purpose is to fit predicted values from the graph embedding. Its output is the districting cost that we want to estimate.
\end{enumerate}
A high-level representation of this architecture is provided in Figure~\ref{fig:gnn}.
It predicts the associated routing cost $\hat{\Phi}(d)$ of a district $d$ given as input and represented as a graph. Detailed information about the three components of the architecture is provided in this section.\\

\noindent
\textbf{Node Embedding Layers.}
Let $G(V,E)$ be the graph representation of the metropolitan area, and $\bm{f}_{vd}$ be the features of each node $v \in V$ for a district $d$.
Formally, a \texttt{GNN} computes a $m$-dimensional features embedding $\bm{\mu}_v \in \mathbb{R}^m$ for each node $v \in V$ in $G$ (i.e., the node embedding). 
The node features $\bm{f}_{vd}$ are aggregated iteratively with the neighboring nodes in the graph. After a predefined number of aggregation steps, the embedding of each node is produced and captures both local and global characteristics of the graph. These operations can be carried out in different ways \citep{dai2016discriminative,velickovic2018graph,hamilton2017inductive}, and public implementations are available for many of these architectures~\citep{wang2019deep,fey2019fast}. Following \cite{dai2017learning}, who considered graph neural networks solving combinatorial problems over graphs, this paper is based on their implementation referred to as \textsc{Structure2Vec} \citep{dai2016discriminative}. We note that other architectures based on GNNs could be used as well.

Let $T$ be the number of aggregation steps, $\bm{\mu}_v^t$ be the node embedding of $v$ obtained after $t$ steps, and $\mathcal{N}(v)$ the set of neighboring nodes of $v \in V$ in $G$. The computation of an embedding $\bm{\mu}_v^{t+1}$ is presented in Equation \eqref{eq:graph_embedding}, where $\bm{\theta}_{1} \in \mathbb{R}^{p \times w}$ 
and $\bm{\theta}_{2} \in \mathbb{R}^{p \times p}$ are tensors of weights that are learned during the training phase, and $\texttt{ReLU}(x) = \max(0,x)$ is a non-linear activation function commonly used in deep neural networks \citep{glorot2011deep}:
\begin{equation}
\label{eq:graph_embedding}
\bm{\mu}_v^{t+1} = \texttt{ReLU}\Big(\bm{\theta}_1 \bm{f}_{vd} + \bm{\theta}_2 \sum_{u \in \mathcal{N}(v)} \bm{\mu}_u^t \Big) \ \ \ \forall t \in \{1,\dots,T\}.
\end{equation}

The idea is to compute the new embedding $\bm{\mu}_v^{t+1}$ as a parametrized sum of the previous embedding ($\bm{\mu}_v^{t}$) and the node features ($\bm{f}_{vd}$).
Following the recommendations of the initial implementation, four aggregation steps are done ($T=4$), each hidden embedding is a vector of $64$ values ($\bm{\mu}_v^{2}, \bm{\mu}_v^{3}, \bm{\mu}_v^{4} \in \mathbb{R}^{64}$). Then, a last aggregation step is performed on the last hidden embedding. This is presented in Equation~\eqref{eq:graph_embedding_out}, where  $\bm{\theta}_{3} \in \mathbb{R}^{k \times p} $ is another weight tensor. In our case, $k$ is set to 1024, which yields a $1024$-dimensional vector as the output embedding ($\bm{\mu}_v^\texttt{out} \in \mathbb{R}^{1024}$):
\begin{equation}
\label{eq:graph_embedding_out}
\bm{\mu}_v^{\texttt{out}} = \texttt{ReLU}\Big(\bm{\theta}_3 \sum_{u \in \mathcal{N}(v)} \bm{\mu}_u^T \Big).
\end{equation}

\noindent
\textbf{Graph Embedding layer.}
Once a vector representation $\bm{\mu}_v^{\texttt{out}}$ has been computed for each vertex~$v$, this information is used to compute $\bm{z}$, a vector representation of the entire graph. This is done by summing together the embedding of each node and applying a non-linear transformation (e.g., \texttt{ReLU}) to the result.
This is also referred to as a \textit{pooling} operation:
\begin{equation}
\label{eq:graph_embedding2}
\bm{z} =  \texttt{ReLU}\Big(  \sum_{v \in V} \bm{\mu}_v^{\texttt{out}} \Big).
\end{equation}

The transformations carried out by Equations (\ref{eq:graph_embedding}--\ref{eq:graph_embedding2}) can be summarized as a parametrized function $\texttt{GNN}:(G \times \mathbb{R}^{8V}) \to \mathbb{R}^{1024}$, which takes as input a graph decorated with eight features at each node, and returns a vector of $1024$ features characterizing the graph.\\

\noindent
\textbf{Fully-Connected Neural Network.}
Finally, the embedding $\bm{z}$ goes through a standard fully connected neural network of two layers, with $100$ neurons for the hidden layer and a single neuron for the output layer. This neural network can be represented as a function $\texttt{FCNN}:\mathbb{R}^{1024} \to \mathbb{R}$, that computes the expected delivery cost $\hat{\Phi}(d)$ inside a district $d$ thanks to the pre-computed graph embedding and
two additional tensor of weights: $\bm{\theta}_4 \in \mathbb{R}^{100 \times 1024}$ and $\bm{\theta}_5 \in \mathbb{R}^{1 \times 100}$.
Assembling all the pieces together, the delivery cost of a district $d$, represented as a graph $G$, and the features
$\bm{f}_{vd}$ for each $v \in V$ are computed as follows:
\begin{equation}
\label{eqprediction}
\hat{\Phi}(d) = \texttt{FCNN}\Bigg(\texttt{GNN}\Big(G, \big[\bm{f}_{vd} ~ \big| ~ v \in V\big]  ~ ; ~  \bm{\theta}_1,\dots, \bm{\theta}_3 \Big) ~ ; ~  \bm{\theta}_4, \bm{\theta}_5 \Bigg).
\end{equation} 

\subsection{Training}
\label{sec:GNN-train}

The network, parameterized by all $\theta$ tensors (171,465 free parameters), is trained via back-propagation using the mean absolute error as loss function (L1 loss). The training is carried out for a maximum of 24 hours or $10,000$ epochs (with a batch size of 64) using Adam optimizer \citep{kingma2014adam}. As an additional early stopping criterion, the training is aborted when no reduction in the loss is observed after $1,000$ consecutive epochs.
We used the default values for the optimizer ($\beta_1=0.9$, $\beta_2=0.999$, and no weight decay), except for the learning rate, which was set to $10^{-4}$ instead of $10^{-3}$ to better stabilize the training. 

\section{Districting Solution Strategy}
\label{sec:districting-algorithm}

The trained \texttt{GNN} provides us with a fast estimation oracle for district costs.
With this, we can now focus our attention on the solution to the districting-and-routing problem. As seen in Section~\ref{sec:problem-statement}, solving the problem through its graph-partitioning mathematical formulation is impractical for most medium- and large-scale instances, due to the difficulty of modeling the objective and the large number of possible districts. To permit experimentation on larger cases, we rely on a fairly standard local search-based approach, which only requires access to the district-cost evaluation oracle. As seen in our numerical experiments and in the Appendix, this method closely matches the results of the exact set-partitioning-based algorithm on small cases solvable to optimality and provides a faster and more scalable solution approach otherwise. Our approach is built upon the Iterated Local Search (ILS) principle \citep{Lourenco2010ILS}, a well-known strategy for guiding local search heuristics beyond local optima. As summarized in Algorithm \ref{alg:ils}, it simply consists of iteratively applying a local search heuristic from an initial solution to attain a local minimum and then perturbing this solution to generate a new starting point for the local search. This process is iterated until a stopping criterion is met. The section briefly describes the different components of the ILS: the construction of an initial solution, the local search moves, and the perturbation scheme.

\begin{figure*}[htbp]
\centering
\begin{minipage}{0.8\textwidth}
\begin{algorithm}[H]
\SingleSpacedXI
\caption{Iterated Local Search (ILS)}
\label{alg:ils}
 $s \gets \text{Generate Initial Solution}$\;
 $s \gets \text{Local Search}(s)$\;
 $s_\textsc{Best}\gets s$\;
 \While{\emph{termination criteria not attained}}{
 $s \gets \text{Pertubation}(s)$\;
 $s \gets \text{Local Search}(s)$\;
 
 \lIf{$\text{\emph{cost}}(s) < \text{\emph{cost}}(s_\textsc{Best})$}{$s_\textsc{Best}\gets s$}
 }
\Return $s_\textsc{Best}$\;
\end{algorithm}
\end{minipage}
\end{figure*}

\subsection{Initial Solution}
\label{sec:initial_solution}

As discussed in Section~\ref{sec:problem-statement}, finding an initial solution of the districting-and-routing problem is NP-hard due to the connectivity and district balance constraints. Consequently, we rely on the solution of an MILP to construct an initial feasible solution, using a compact network-flow formulation of the balanced connected $k$-partition problem from \citet{miyazawa2020flow}.

Given the initial graph $G(V, E)$ representing the adjacency between the different BUs, the network flow formulation uses an auxiliary directed graph $H(V \cup S, A)$ in which $S$ is a set of~$k$ source vertices, one for each district. The set of arcs $A$ is obtained by including two opposite arcs $(u,v)$ and $(v,u)$ for each edge $\{u,v\} \in E$, as well as an arc between each $u \in S$ and $v \in V$. Binary design variables $y_{ij}$ take the value $1$ when $(i,j)$ carries some flow, and $f_{ij}$ is the amount of flow. With this, the feasibility model seeks a solution that satisfies:
\begin{align}
 \label{flow-equation:flow-balance-constraint}  & \hspace*{-0.2cm} \sum_{u \in V \cup S} f_{uv} - \sum_{u \in V \cup S} f_{vu} = 1 & & \forall v \in V & \\
 \label{flow-equation:flow-upper-bound} & f_{uv} \leq n_\textsc{u} y_{uv} & & \forall u,v \in V \cup S & \\
 \label{flow-equation:source-min-flow} & \sum_{v \in V}f_{sv} \geq n_\textsc{l} & & \forall s \in S & \\
 \label{flow-equation:symmetry-constraint} & \sum_{v \in V} f_{sv} \leq \sum_{v \in V} f_{s+1,v} & & \forall s \in \{0, \dots, k-1\} & \\
 \label{flow-equation:source-minimum} & \sum_{v \in V} y_{sv} \leq 1 & & \forall s \in S & \\
 \label{flow-equation:edges-coming-maximum} & \hspace*{-0.2cm} \sum_{u \in V \cup S} y_{uv} \leq 1 & & \forall v \in V & \\
 \label{flow-equation:y-as-bool} & y_{uv} \in \{0,1\} & & \forall u,v \in V \cup S & \\
 \label{flow-equation:flow-as-continuous} & f_{uv} \in \mathbb{R}^{+} & & \forall u,v \in V \cup S.
\end{align}

This model is constructed so that sources $s \in S$ representing districts send flow towards vertices $v \in V$ representing the BUs. Constraints (\ref{flow-equation:flow-balance-constraint}) ensure flow conservation and guarantee that each BU vertex $v \in V$ receives one unit of flow. Constraint (\ref{flow-equation:flow-upper-bound}) sets a limit on the maximum flow of an arc and consequently enforces the upper limit $n_\textsc{u}$ on the number of BUs per district. Constraint (\ref{flow-equation:source-min-flow}) ensures that each source vertex feeds flow to at least $n_\textsc{l}$ BUs. Constraint~(\ref{flow-equation:symmetry-constraint}) orders the sources by flow amount to eliminate symmetry. Finally, Constraints~(\ref{flow-equation:source-minimum}--\ref{flow-equation:edges-coming-maximum}) uniquely match each source vertex of $S$ with a single receiving BU vertex of $V$ (which will then forward flow to other connected BU vertices -- ensuring the connectivity constraints) and ensure that each BU is connected to a single BU or source. 
Solving this MILP on graphs with a few hundred BUs (as in our cases) takes seconds using modern solvers (Gurobi v12.0.0). Then, any feasible solution of this model can be converted into an initial solution of the districting-and-routing problem by assigning each BU vertex $v \in V$ to the (unique) district vertex $s \in S$ from which its flow originated.

\subsection{Local Search}

After the initial solution construction and after each perturbation, our solution approach applies a local search (LS) procedure. As summarized in Algorithm \ref{alg:local-search}, this LS follows a first-improvement policy that consists of exploring a neighborhood in a random order and directly applying any improving move until a local minimum is attained. For any solution $s$, the neighborhood $\mathcal{N}(s)$ is defined as the set of solutions that can be attained by applying any of these two moves: (1) \textsc{Relocate}$(u,d)$, which reassigns a BU $u$ to a new district $d$, and (2) \textsc{Swap}$(u,v)$, which exchanges the districts of two BUs $u$ and $v$. Therefore, the neighborhood size is such that $|\mathcal{N}(s)| = \mathcal{O}(kn + n^2) = \mathcal{O}(n^2)$, but moves leading to infeasible solutions w.r.t. district-cardinality and connectivity constraints can be directly disregarded. This permits several speedup strategies. Firstly, moves can be limited to pairs of districts $d_i$ and $d_j$ that are connected, such that there exist two neighbors BUs $(u,v) \in E$ with $u \in d_i$ and $v \in d_j$. Besides this, we can define border $\mathcal{B}_{ij}$ as the set of BUs of a district $d_i$ that have at least one neighbor BU belonging to district $d_j$. We can then further limit the set of moves between districts $d_i$ and $d_j$ to BUs $u \in \mathcal{B}_{ij}$ and $v \in \mathcal{B}_{ji}$.

\begin{figure*}[htbp]
\centering
\begin{minipage}{0.9\textwidth}
\begin{algorithm}[H]
\SingleSpacedXI
\caption{Local search procedure}\label{alg:local-search}

\While{\emph{Improvement found}}
{
 \For{\emph{all pairs of \emph{connected districts} $
 \{d_i, d_j\}$ in random order}}
 {
 Find and apply, if improving, the best move on $s$ among:
 
 \hspace*{1em} $\bullet$ all feasible \textsc{relocate} moves for $\{(u, d_j) \, | \, u \in \mathcal{B}_{ij}\} \cup \{(v, d_i) \, | \, v \in \mathcal{B}_{ji}\}$; 
 
 \hspace*{1em} $\bullet$ all feasible \textsc{swap} moves for $ \{(u,v) \, | \, u \in \mathcal{B}_{ij}, \, v \in \mathcal{B}_{ji}\}$; 
 }
}
\end{algorithm}
\end{minipage}
\end{figure*}

\begin{figure}%
 \centering
 \includegraphics[width=0.85\textwidth]{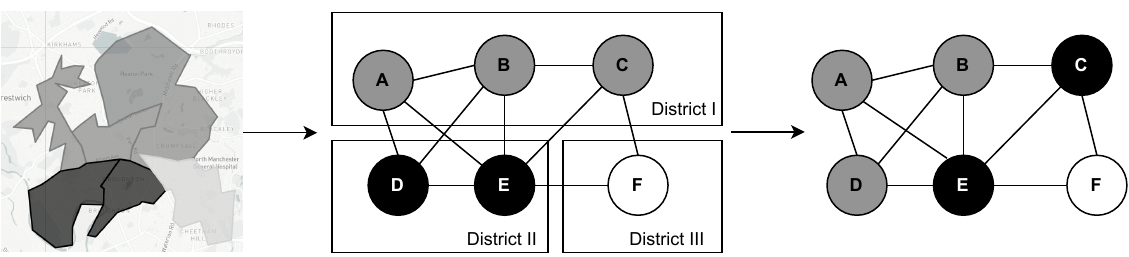}
 \caption{Local Search example}
 \label{fig:local-search-example}%
\end{figure}

Figure \ref{fig:local-search-example} illustrates this local search procedure. It represents a city with six BUs that must be partitioned into three districts, along with the associated graph $G$. For the current solution, the possible moves between districts 1 and 2 include \textsc{Relocate} moves $(A,II), (B,II), (C,II), (D,I), (E,I)$ and \textsc{Swap} moves $(A, D), (A, E), (B, D), (B, E), (C, D), (C, E)$. For example, applying \textsc{Swap}(C, D) leads to the rightmost solution. Though $C$ and $D$ are not adjacent, it is possible to swap them since  $\text{C} \in \mathcal{B}_{\text{I,II}}$ and $\text{D} \in \mathcal{B}_{\text{II,I}}$.

\subsection{Perturbation}

The perturbation operator is an important building block of any ILS. It must permit escaping from local optima without losing too much structural information on good solutions. Our perturbation procedure simply uses the same set of moves as the local search but applies any feasible move with a probability $p_{\textsc{rm}}$ regardless of its cost. Parameter $p_{\textsc{rm}}$ 
drives the amplitude of the perturbation. Its calibration is discussed in the Appendix.

\section{Experimental Analyses}
\label{sec:experiments}

For our experimental analyses, we focus on districting-and-routing problem instances that occur in five metropolitan areas in the UK (London, Bristol, Manchester, Leeds, and the West Midlands) with very diverse geographical characteristics. The goal of our experiments is twofold. First, we aim to evaluate the proposed GNN's routing-cost estimation accuracy. Next, we measure the extent to which cost-estimation accuracy impacts the ability to make good strategic decisions, i.e., to what extent solutions of districting problems using the GNN or other cost predictors differ in terms of their partitions and operational efficiency.

In the rest of this section, we discuss the data collection and generation of test instances (Section~\ref{subsec:dataset}), as well as the baseline methods considered for routing-cost estimation (Section~\ref{subsec:baselines-routing-cost}). Then, we analyze the accuracy of the different prediction cost-estimation methods (Section~\ref{subsec:results-predictive-methods}), and measure the impact of different estimation algorithms when optimizing districts as well as the characteristics (e.g., compactness) of the partitions thereby generated (Section~\ref{subsec:impact-on-districting}).

\subsection{Data Collection and Experimental Setup}
\label{subsec:dataset}

\noindent
\textbf{Geographical data and test instances.}
We base our studies on five metropolitan areas located in the UK (London, Bristol, Manchester, Leeds, and the West Midlands). The selection of these areas was driven by their diversity and availability of precise geographical boundaries from \url{https://movement.uber.com} as well as population statistics from the UK government 2018 census database \citep{CensusData2018}. The BUs correspond to Middle Layer Super Output Areas (MSOAs), which are designed to contain roughly the same population (each MSOA contains over 5,000 inhabitants, and 8,000 inhabitants on average over the UK). Finally, to obtain data sets with a different number of BUs ($n = \{60, 90 ,120\}$), we selected a center point in each region and retained the $n$ closest BUs.

Table~\ref{tab:cities-stats} provides general statistics (population, area, density, and compactness) on the $n=120$ BUs of each metropolitan area. The area of each BU was measured using the Monte Carlo method with 50,000 samples. Compactness scores have been obtained using Reock's formula \citep{Young_1988}, which divides the area of the BU by the area of the smallest circumscribed circle. Therefore, this compactness measure assigns higher scores to areas that are more circularly shaped.

\begin{table}[!htbp]
\setlength\tabcolsep{10pt}
\renewcommand{\arraystretch}{1.2}
\centering
\caption{Population and geography statistics for the BUs}\label{tab:cities-stats}
\scalebox{0.75}
{
\begin{tabular}{lll|M{1.8cm}M{1.8cm}M{1.8cm}cc}
\toprule
& & & \textbf{Bristol}& \textbf{Leeds}& \textbf{London}& \textbf{Manchester}& \textbf{West-Mid.}\\
\hline
\multirow{5}{*}{\textbf{Population (thousands)}}
 & & average& 8.32 & 7.53 & 9.68 & 8.69 & 8.71\\
 & & std& 2.29 & 1.64 & 2.03 & 2.36 & 2.06\\
 & & min& 5.55 & 5.20 & 6.58 & 5.26 & 5.44\\
 & & median& 7.68 & 7.26 & 9.40 & 8.32 & 8.19\\
 & & max& 18.16 & 14.06 & 16.17 & 15.87 & 17.12\\
\hline
\multirow{5}{*}{\textbf{Area (km$^2$)}}
 & & average& 10.34 & 4.67 & 0.75 & 2.22 & 1.81\\
 & & std& 22.42 & 7.22 & 0.47 & 1.17 & 0.75\\
 & & min& 0.63 & 0.35 & 0.30 & 0.59 & 0.53\\
 & & median& 1.99 & 2.51 & 0.64 & 1.90 & 1.68\\
 & & max& 171.21 & 51.79 & 3.58 & 6.69 & 4.39\\
\hline
\multirow{5}{*}{\textbf{Density (thousands/km$^2$)}}
 & & average& 3.92 & 3.61 & 15.15 & 4.87 & 5.51\\
 & & std& 2.89 & 3.41 & 4.91 & 2.58 & 2.41\\
 & & min& 0.06 & 0.14 & 2.76 & 1.12 & 1.99\\
 & & median& 3.74 & 2.90 & 15.17 & 4.41 & 5.21\\
 & & max& 12.72 & 25.20 & 28.27 & 16.36 & 17.06\\
\hline
\multirow{5}{*}{\textbf{Compactness}}
 & & average& 0.40 & 0.43 & 0.43 & 0.42 & 0.44\\
 & & std& 0.10 & 0.09 & 0.09 & 0.10 & 0.09\\
 & & min& 0.19 & 0.25 & 0.21 & 0.20 & 0.22\\
 & & median& 0.40 & 0.43 & 0.43 & 0.42 & 0.44\\
 & & max& 0.65 & 0.67 & 0.63 & 0.66 & 0.61\\
\bottomrule
\end{tabular}
}
\end{table}

As seen in Table~\ref{tab:cities-stats}, there are important differences in geography and population between the BUs of the five metropolitan areas. Bristol and Leeds BUs generally have smaller population densities and larger geographical extensions, whereas London has the densest and smallest BUs. BUs areas are also very disparate, ranging from 0.30 km$^2$ (in London) to 171.21 km$^2$ (in Bristol). Finally, compactness measures appear to be roughly similar for the five cities.

Besides the characteristics of the BUs, other factors impact the structure of the district-and-routing problems: the admissible range for the number of BUs in each district [$n_\textsc{l},n_\textsc{u}]$, and the location of the depot. For $t \in \{3,6,12,20,30\}$, we generated five different configurations for the district-size constraints by setting $n_\textsc{l} = \lfloor 0.8 \times t \rfloor$ and $n_\textsc{u} = \lceil 1.2 \times t \rceil$, and setting the number of districts to $k = \lfloor n / t \rfloor$. Then, for each configuration, the probability $\kappa$ that an inhabitant will make a request has been set to $\kappa = 96/ (8000t)$ where $96$ corresponds to the targeted number of requests in the related routing problem, and $8000$ corresponds to the average population of a BU. The values for parameter $\kappa$ have been selected to reflect realistic scenarios with delivery routes that cover close to a hundred stops (typical for parcel deliveries, as seen in \citealt{Holland2017}). Finally, five possible locations for the depot were considered: at the center of the metropolitan area (\texttt{C}), at the north-east (\texttt{NE}), at the north-west (\texttt{NW}), at the south-east (\texttt{SE}), and at the south-west (\texttt{SW}) of the center ($D = \{\texttt{C},\texttt{NE},\texttt{NW},\texttt{SE},\texttt{SW}\}$).
All these factors are summarized in Table \ref{table:config}. For each of the five considered metropolitan areas, we therefore generated $3 \times 5 \times 5 = 75$ instances covering all the possible combinations of these factors.\\

\begin{table}[!htbp]
\setlength\tabcolsep{10pt}
\renewcommand{\arraystretch}{1.2}
\centering
\caption{Summary of the different instance parameters}
\scalebox{0.9}
{
\begin{tabular}{ll}
\hline
Factor & Values \\ 
\hline 
Total number of BUs & $n \in \{60,90,120\}$ \\
Target number of BUs in a district & $t \in \{3,6,12,20,30\}$ \\
Depot location & $\{\texttt{C},\texttt{NE},\texttt{NW},\texttt{SE},\texttt{SW}\}$ \\
Minimum number of BUs in a district & $n_\textsc{l} = \lfloor 0.8 \times t \rfloor$ \\
Maximum number of BUs in a district & $n_\textsc{u} =\lceil 1.2 \times t \rceil$ \\
Number of districts & $k = \lfloor n / t \rfloor$ \\
Request probability & $\kappa = 96/(8000t)$ \\ 
\hline
\end{tabular}
}
\label{table:config}
\end{table}

\noindent
\textbf{Scenarios and solution evaluations.}
The districting-and-routing problem has an objective function that must be measured over a large number of scenarios. Consequently, to standardize the evaluation procedures among the different evaluation and solution methods, it was essential to (i) create separate demand scenarios (positions of customers' demands) for training and testing and (ii) adopt a standard for solution measurement. We did this by sampling for each BU $S_\textsc{Train} = 500$ and $S_\textsc{Test} = 500$ random scenarios for training and testing, respectively. Let $X^\textsc{T}_{it}$ be the set of customer requests, characterized by their positions, for any BU $i$ and scenario $s$ for $\textsc{T} \in \{\textsc{Train},\textsc{Test}\}$. All methods have access only to the training set of scenarios during the learning and solution process. We reserve the other scenarios for the final evaluation and comparison of the solutions. For each demand scenario of the BUs, we obtain a corresponding demand scenario for any district by compounding all the demand locations of the BUs it contains. With this, the training (and respectively testing) cost of a solution $\mathcal{S}$ described as a set of districts, is defined as the average cost of the TSP tours over the scenarios, calculated as:
\begin{equation}
\Phi^\textsc{T}_\textsc{saa}(\mathcal{S}) = \frac{1}{S_\textsc{T}} \sum_{s=1}^{S_\textsc{t}} \sum_{d \in \mathcal{S}} C_\textsc{tsp}\left(\bigcup_{i \in d} X^\textsc{T}_{is} \right) ,
\end{equation}
for $\textsc{T} \in \{\textsc{Train},\textsc{Test}\}$, and where $C_\textsc{tsp}(X)$ calculates the cost of a TSP tour visiting the depot and all customers in $X$. We rely on the Lin-Kernighan (LKH) algorithm (from \citealt{Helsgaun2000}, available at \url{http://webhotel4.ruc.dk/~keld/research/LKH/}) to measure TSP distances. The reason for this choice is that this algorithm is very fast (0.1 seconds per run) and finds the optimal solution almost systematically in our cases with roughly 100 customer requests, so there is no need to employ time-consuming exact approaches that would limit our experimentation capabilities.\\

\noindent
\textbf{Computational Environment.} 
The calculation of TSP scenarios to estimate district costs has been conducted on Intel Gold 6148 Skylake 2.4 GHz CPUs with 80GB RAM. The rest of the 
experimental analyses (including training, calibration, and optimization) have been conducted on a computer with an Intel E5-2683 v4 Broadwell 2.1GHz CPU, 124GB RAM, and an Nvidia P100 Pascal (12GB HBM2 memory) GPU. To implement the \texttt{GNN}, we use the official \texttt{structure2vec} implementation via its \texttt{Python} interface, available at \url{https://github.com/Hanjun-Dai/pytorch_structure2vec}.
The ILS is implemented in \textsc{C++}, compiled with \textsc{G++ v9.3.0}.
Finally, to ensure experimental reproducibility, all data, scripts, and source code needed to reproduce our experiments are available in an open-source repository at \url{https://github.com/vidalt/Districting-Routing}.

\subsection{Estimation Methods for Routing Costs}
\label{subsec:baselines-routing-cost}

We identified three main estimation approaches for routing costs in districting problems: \citeauthor{Beardwood1959}'s formula extended by \citeauthor{daganzo1984distance} (called BHHD in the rest of this paper), the variant of this formula by \citeauthor{Figliozzi2007} (FIG), and the shallow neural network designed by \citeauthor{Kwon1995} (SNN). We will rely on these methods for our experimental comparisons.\\

\noindent
\begin{itemize} 
\item\textbf{BHHD}:
Let $d$ be a district composed of $n$ BUs. For any BU $i$, we recall that $a_i$ is the area of the BU and $\xi_i$ is its population. Moreover, $\kappa$ is the probability of a delivery request for any inhabitant. With this, the total area of the district is calculated as $A_d = \sum_{i \in d} a_i$, and the expected number of deliveries within the district is $R_d = \kappa \sum_{i \in d} \xi_i$. Extending the BHH formula to account for the distance from the depot, as in \cite{daganzo1984distance} gives:
\begin{equation}
\Phi_{\textsc{bd}}(d) = \alpha \sqrt{A_d R_d} + 2 \Delta_d,
\label{eq:bd-daganzo-loos}
\end{equation}
where $\alpha$ is a hyper-parameter that needs calibration, and $\Delta_d$ is the average distance between the depot and a request, calculated by Monte-Carlo estimation on the training scenario set. In our experiments, we set $\alpha$ to minimize the mean-squared error over a training set (discussed at the end of this section) through least-squares regression. 
We note that the BHH formula converges almost surely to the true expected distance in an asymptotic regime where the number of delivery requests tends towards infinity. As a consequence, approximation errors naturally occur as we deviate from this assumption.\\

\item\textbf{FIG}:
This continuous approximation formula of \cite{Figliozzi2007} is a direct extension of~Equation~\eqref{eq:bd-daganzo-loos}, which was designed to cope with the practical non-uniformity of real observed demands. The formula is defined in Equation~\eqref{eq:fig-loss}, and includes four hyper-parameters ($\alpha_1$, $\alpha_2$, $\alpha_3$, $\alpha_4$) requiring calibration. As previously, $\Delta_d$ represents the average distance between the depot and the customer requests. As previously, the values of the hyperparameters are selected through least-squares regression.
\begin{equation}
\texttt{FIG}(d) = \alpha_1 \sqrt{A_d R_d} + \alpha_2 \Delta_d + \alpha_3 \sqrt{\frac{A_d}{R_d}} + \alpha_4\\
\label{eq:fig-loss}
\end{equation}

\item\textbf{SNN}:
This estimator is based on a neural network with one hidden layer of three neurons. Five features are used as input to the network:
\begin{enumerate}[(1)]
 \item The expected number of deliveries in the district ($R_d$);
 \item The ratio between the length $d^l$ and the height $d^h$ of a minimal-area rectangle covering the district: ($d^l/d^h$); 
 \item The average distance of a customer from the district $d$ to the depot ($\Delta_d$);
 \item Feature (2) divided by Feature (1).
 \item The BHH distance estimate ($\sqrt{A_d R_d}$).
\end{enumerate}

Let $\bm{x}_d \in \mathbb{R}^5$ be the vector of these five features for a district $d$. The routing cost estimation is given by Equation \eqref{eq:kwon-nn} where $\bm{w}_1 \in \mathbb{R}^{3 \times 5}$ and $\bm{w}_2 \in \mathbb{R}^{1 \times 3}$ are the network weights, and where $\bm{b}_1 \in \mathbb{R}^{3}$ and $b_2 \in \mathbb{R}$ are the biases. This leads to a total of $22$ free parameters to learn.
Finally, $\gamma : \mathbb{R} \to [0,1]$ is a non-linear activation function.
\begin{equation}
\label{eq:kwon-nn}
\texttt{SNN}(d) = \bm{w}_2 \gamma\Big(\bm{w}_1 \bm{x}_d + \bm{b}_1 \Big) + b_2
\end{equation}

\cite{Kwon1995} used a sigmoid activation function and standard backpropagation for training.
However, this early architecture did not benefit from the subsequent advances in deep learning and performed poorly compared to the other baselines. To ensure a fair comparison under modern training conditions, we replaced the sigmoid activations with \texttt{ReLU} functions and trained the network using the Adam optimizer \citep{kingma2014adam} with a learning rate of $10^{-3}$. The SNN model is therefore trained for $50,000$ epochs, considering the mean squared error as the loss function.\\
\end{itemize}

\noindent
\textbf{Calibration and Training.}
To calibrate and train the different cost-estimation methods, we sampled for each instance 9,000 random connected districts respecting the size constraints $[n_\textsc{l},n_\textsc{u}]$.
For each of these districts, we calculated the expected TSP cost by SAA on the \textsc{Train} demand scenarios. For GNN and SNN, we further subdivided this set into 8,000 districts for training and 1,000 districts for validation to control convergence. Training the GNN takes two hours on average, i.e., a substantial but reasonable computational effort since strategic districting decisions hold for a long time. Afterwards, district cost inference takes on average 8 milliseconds on a single CPU thread.

\subsection{Analysis -- Prediction Accuracy}
\label{subsec:results-predictive-methods}

Our first set of experiments aims to evaluate the accuracy of the different models (\texttt{BHHD}, \texttt{FIG}, \texttt{SNN}, as well as the proposed \texttt{GNN}) for estimating the routing costs. Therefore, we use the trained models described in the previous section and evaluate them on an additional set of 1,000 random districts that are distinct from the ones used during training. We compare the estimates with ``ground truth'' values of the routing costs obtained again by SAA over the $500$ \textsc{Test} demand scenarios. In the rest of this section, we analyze how the proposed estimation approaches deviate from the ground-truth measurements in terms of their Root-Mean-Square Error (RMSE) \myblue{and Mean Absolute Percentage Error (MAPE)}. For a given test instance, these values are calculated as:
\begin{equation}
\text{RMSE} = \sqrt{\frac{1}{|\mathcal{D}|}
\sum_{d \in \mathcal{D}} \left(\Phi(d) - \Phi^\textsc{test}_\textsc{saa}(d)\right)^2}, 
\hspace*{1cm}  \myblue{\text{MAPE} = \frac{100}{| \mathcal{D}|} \sum_{d \in \mathcal{D}}
    \left| \frac{\Phi(d) - \Phi^\text{TEST}_{\text{SAA}}(d)}
    {\Phi^\text{TEST}_{\text{SAA}}(d)} \right|,}
\end{equation}
where $\mathcal{D}$ is the set of 1,000 evaluation districts for this instance, and $\Phi(d)$ is the cost estimation provided by the considered method on a district $d$.

Table~\ref{tab: cost-estimators-mse} compares the RMSE of the different cost-estimation approaches. Each line corresponds to the results for the instances with a certain number of BUs ($n$) and district-size target ($t$), averaging the RMSE over the five different metropolitan areas and five possible depot-position configurations. The columns provide the characteristics of the instances, the average cost of the districts evaluated through SAA, and finally, the RMSE of the different cost estimation approaches.
Then, Table~\ref{tab: cost-estimators-mse-2} focuses specifically on the medium-scale case of $n=90$ and $t=12$, with additional detailed results for each metropolitan area and depot configuration.

\begin{table}[!htbp]
    \setlength\tabcolsep{12.5pt}
    \renewcommand{\arraystretch}{1.1}
    \centering
    \caption{Accuracy of the different estimation approaches}\label{tab: cost-estimators-mse}
    \scalebox{0.85}
    {
    \begin{tabular}{ccccccccccc}

    \toprule
    \multirow{2}{*}{\textbf{$\bm{n}$}} & \multirow{2}{*}{$\bm{t}$} & \multirow{2}{*}{$\hat{\Phi}^\textsc{test}_\textsc{saa}$} & \multicolumn{4}{c}{\textbf{RMSE}} & \multicolumn{4}{c}{\textbf{\myblue{MAPE}}} \\
     \cmidrule(rl){4-7}
     \cmidrule(rl){8-11}
    & & & \texttt{GNN} & \texttt{BHHD} & \texttt{FIG} & \texttt{SNN} & \texttt{GNN} & \texttt{BHHD} & \texttt{FIG} & \texttt{SNN}\\

    \midrule
    \multirow{5}{*}{60}

    &3 &43.16 $\pm$ 0.48 &\textbf{1.48} &2.27 &2.09 &2.05  & \textbf{2.75} & 3.81 & 3.60 & 3.60 \\
    &6 &48.21 $\pm$ 0.66 &\textbf{1.39} &3.43 &3.24 &3.02  & \textbf{2.20} & 5.13 & 4.87 & 4.91 \\
    &12 &52.06 $\pm$ 0.80 &\textbf{2.10} &4.36 &4.15 &4.00 & \textbf{2.85} & 6.19 & 5.78 & 5.79  \\
    &20 &56.42 $\pm$ 0.95 &\textbf{2.39} &5.18 &4.89 &4.73 & \textbf{2.83} & 6.62 & 6.11 & 6.18  \\
    &30 &61.14 $\pm$ 1.04 &\textbf{2.08} &4.73 &4.50 &4.31 & \textbf{2.41} & 5.92 & 5.49 & 5.56  \\

    \midrule
    \multirow{5}{*}{90}

    & 3 &52.43 $\pm$ 0.53 &\textbf{1.56}  &2.70 &2.49 &2.47 &\textbf{2.47} &3.68 &3.53 &3.55 \\
    & 6 &57.03 $\pm$ 0.69  &\textbf{1.88} &4.00 &3.75 &3.60 &\textbf{2.52} &4.93 &4.66 &4.73 \\
    & 12 &61.19 $\pm$ 0.83 &\textbf{2.62} &5.12 &4.80 &4.60 &\textbf{3.05} &6.11 &5.60 &5.76 \\
    & 20 &66.67 $\pm$ 1.02 &\textbf{3.00} &6.17 &5.74 &5.50 &\textbf{3.00} &6.75 &6.12 &6.27 \\
    & 30 &72.63 $\pm$ 1.14 &\textbf{2.84} &6.15 &5.78 &5.50 &\textbf{2.75} &6.33 &5.91 &6.00 \\

    \midrule

    \multirow{5}{*}{120}

    & 3 &61.61 $\pm$ 0.60  &\textbf{1.68} &3.27 &3.15 &3.10 & \textbf{2.27} &3.60 &3.46 &3.49 \\
    & 6 &63.25 $\pm$ 0.69  &\textbf{2.09} &4.04 &3.88 &3.76 & \textbf{2.55} &4.57 &4.35 &4.42 \\
    & 12 &67.57 $\pm$ 0.84 &\textbf{2.83} &5.28 &5.00 &4.84 & \textbf{3.06} &5.71 &5.29 &5.38 \\
    & 20 &74.44 $\pm$ 1.07 &\textbf{3.34} &6.36 &6.03 &5.90 & \textbf{3.15} &6.23 &5.77 &5.84 \\
    & 30 &80.81 $\pm$ 1.21 &\textbf{3.39} &6.82 &6.52 &6.21 & \textbf{2.91} &6.25 &5.78 &5.88 \\

    \midrule
    \multicolumn{2}{c}{Average} &59.84 $\pm$ 0.81 &\textbf{2.31} &4.66  &4.40 &4.24 &\textbf{2.72} &5.46 &5.09 &5.16\\
    \bottomrule
    \end{tabular}
    }
    \end{table}

\begin{table}[!htbp]
    \setlength\tabcolsep{8pt}
    \renewcommand{\arraystretch}{1.1}
    \caption{Impact of the depot location and metropolitan area, for $\bm{n=90}$ and $\bm{t=12}$} \label{tab: cost-estimators-mse-2}
        \hspace*{-0.5cm}
    \scalebox{0.85}
    {
    \begin{tabular}{ccccccccccc}
    \toprule
    \multirow{2}{*}{\textbf{Depot}} & \multirow{2}{*}{\textbf{Metropolitan Area}} & \multirow{2}{*}{$\hat{\Phi}^\textsc{test}_\textsc{saa}$} & \multicolumn{4}{c}{\textbf{RMSE}} & \multicolumn{4}{c}{\textbf{\myblue{MAPE}}} \\
         \cmidrule(rl){4-7}
     \cmidrule(rl){8-11}
    & & & \texttt{GNN} & \texttt{BHHD} & \texttt{FIG} & \texttt{SNN} & \texttt{GNN} & \texttt{BHHD} & \texttt{FIG} & \texttt{SNN}\\
    \midrule
    \multirow{5}{*}{\shortstack{Central Depot \\ $\{\texttt{C}$\}}}

    &\texttt{Bristol} &67.14 $\pm$ 1.19 &\textbf{4.04} &8.38 &7.99 &6.56  & \textbf{4.38} & 8.74 & 9.46 & 7.97 \\
    &\texttt{Leeds} &48.71 $\pm$ 0.97 &\textbf{3.14} &5.75 &5.25 &5.29  & \textbf{4.32} & 8.14 & 8.12 & 7.35 \\
    &\texttt{London} &28.58 $\pm$ 0.59 &\textbf{1.75} &3.25 &3.19 &3.11  & \textbf{4.02} & 8.43 & 7.90 & 8.25 \\
    &\texttt{Manchester} &43.13 $\pm$ 0.89 &\textbf{2.70} &4.51 &4.43 &4.45  & \textbf{4.45} & 8.05 & 7.35 & 7.43 \\
    &\texttt{West-Midlands} &40.59 $\pm$ 0.84 &\textbf{2.24} &4.06 &3.90 &3.88  & \textbf{3.73} & 8.17 & 7.21 & 7.06 \\

    \midrule
    \multirow{5}{*}{\shortstack{Off-Centered \\ Depot \\ $\{\texttt{NE},\texttt{NW},\texttt{SE},\texttt{SW}\}$}}

    &\texttt{Bristol} &105.59 $\pm$ 1.13     &\textbf{3.72} &8.96 &8.06 &7.01  & \textbf{2.25} & 5.93 & 5.30 & 4.50\\
    &\texttt{Leeds} &71.78 $\pm$ 0.86        &\textbf{2.67} &5.70 &5.09 &5.26  & \textbf{2.47} & 5.28 & 5.13 & 4.97\\
    &\texttt{London} &37.03 $\pm$ 0.54       &\textbf{1.69} &2.98 &2.93 &2.90  & \textbf{3.12} & 5.96 & 5.61 & 5.78\\
    &\texttt{Manchester} &57.60 $\pm$ 0.80   &\textbf{2.56} &4.14 &4.09 &4.11  & \textbf{3.05} & 5.41 & 5.06 & 5.19\\
    &\texttt{West-Midlands} &53.39 $\pm$0.76 &\textbf{2.30} &3.70 &3.66 &3.62  & \textbf{2.96} & 5.19 & 4.91 & 5.04\\

    \bottomrule
    \end{tabular}
    }
    \end{table}

\myblue{
As seen in these tables, with an overall RMSE of 2.31 (respectively, a MAPE of 2.72), the quality of the estimates obtained with the proposed \texttt{GNN} is vastly superior to that of all the other methods: \texttt{BHHD} with a RMSE of 4.66 (MAPE 5.46), \texttt{FIG} with a RMSE of 4.40 (MAPE 5.09), and \texttt{SNN} with a RMSE of 4.24 (MAPE 5.16).} The error committed also appears to increase with the number of BUs in the metropolitan area~($n$) as well as the target number of BUs~($t$) in each district. This is due to two factors. Firstly, the magnitude of the predicted values (long-term operational cost of the districts) grows with $t$, and also with $n$ to a lesser extent (instances with larger $n$ include less-populated BUs that are located farther away from city centers). Therefore, as these values grow larger, more estimation error is generally committed. Secondly, the number of possible districts grows exponentially in $n$ and $t$, such that the universe of possible inputs quickly grows, and it becomes harder to learn the target.

Considering the results of Table~\ref{tab: cost-estimators-mse-2}, we again notice that \texttt{GNN} estimations are far more accurate than the other approaches for all metropolitan areas and depot configurations (close or away from the city center). The magnitude of the RMSE depends on the area. Indeed, metropolitan areas such as Bristol have a lower population density and longer tours within the districts, leading to generally higher operational-cost values. In this situation, it is natural for the error to grow with the magnitude of the values reported in column $\hat{\Phi}^\textsc{test}_\textsc{saa}$.

\begin{figure}[htbp]
\setlength\tabcolsep{10pt}
\renewcommand{\arraystretch}{1.2}
\centering
\scalebox{0.78}
{
\begin{tabular}{m{0.165\columnwidth} m{0.425\columnwidth} m{0.425\columnwidth}}
& \multicolumn{1}{c}{\hspace*{0.5cm}\textbf{Bristol}} & \multicolumn{1}{c}{\hspace*{0.5cm}\textbf{London}} \\

($n=90$, $t=3$) & \includegraphics[width=1.\linewidth]{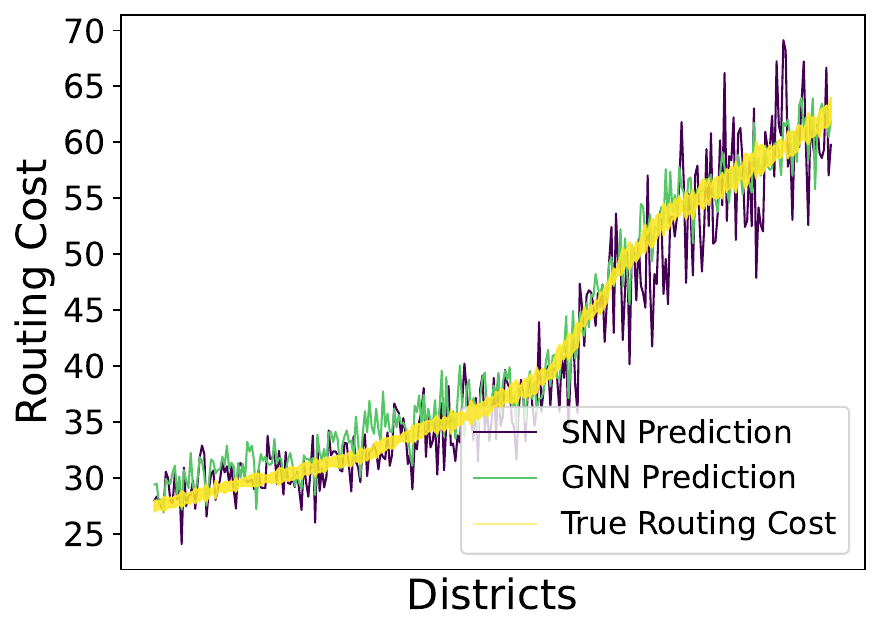} & \includegraphics[width=1.\linewidth]{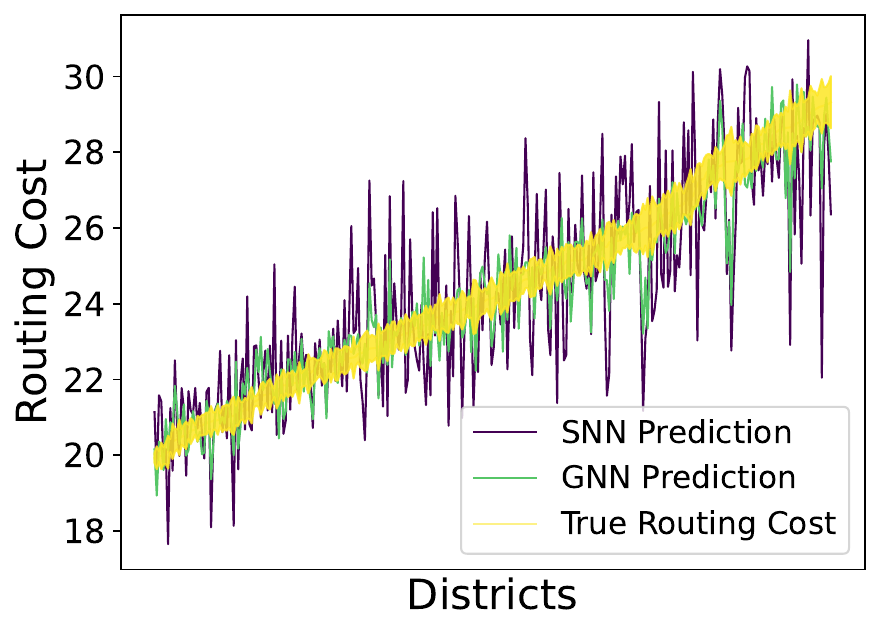} \\
($n=90$, $t=12$) & \includegraphics[width=1.\linewidth]{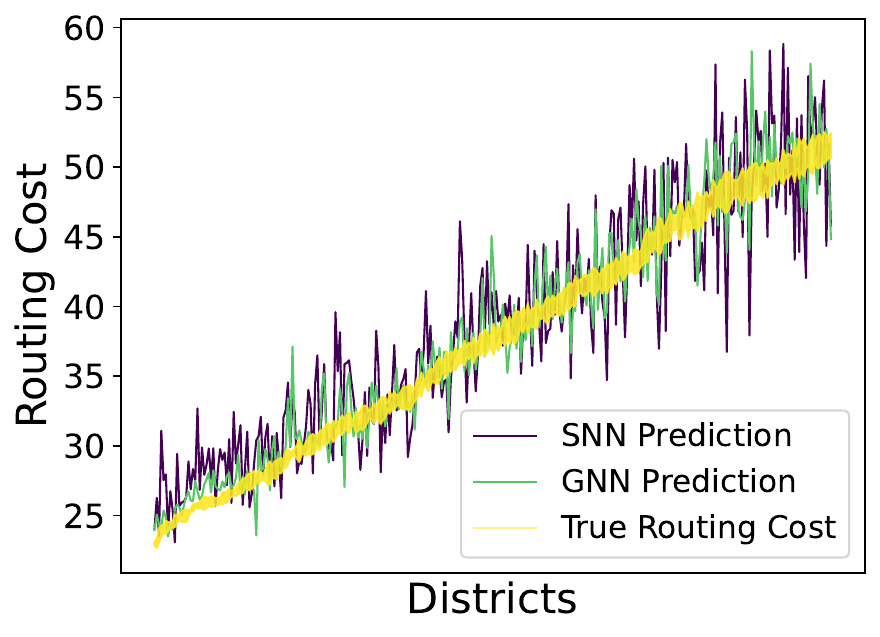} & \includegraphics[width=1.\linewidth]{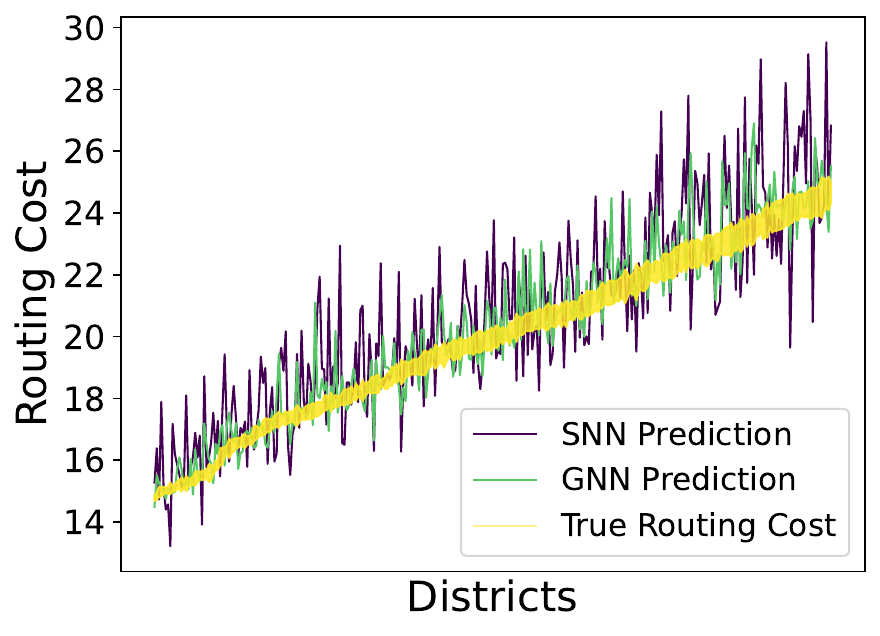} \\
($n=120$, $t=12$) & \includegraphics[width=1.\linewidth]{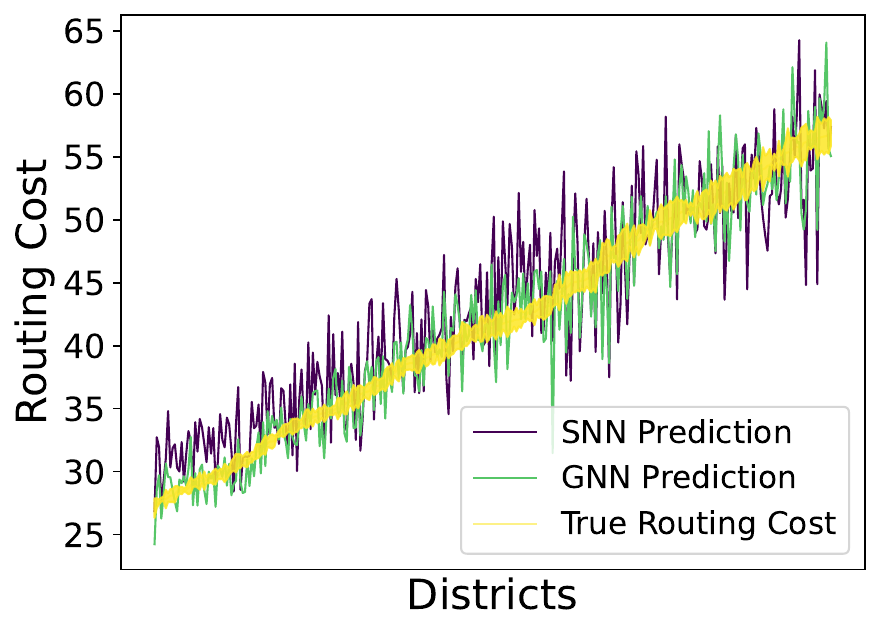} & \includegraphics[width=1.\linewidth]{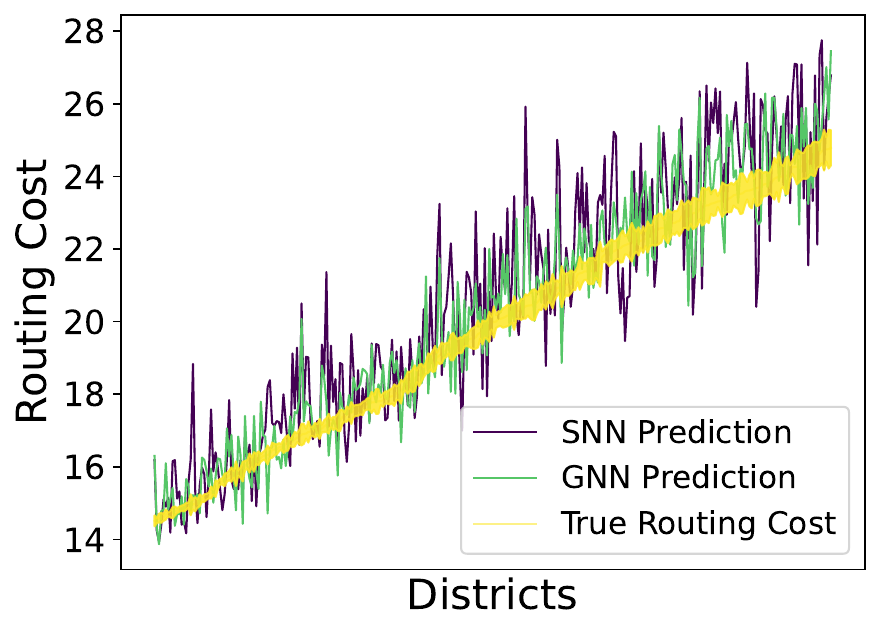} \\
\end{tabular}
}
\caption{Visualization of the ground truth and estimated district costs on a subset of districts for Bristol and London (with a central depot).} \label{fig:snnxgnn-mse}
\end{figure}

Finally, Figure~\ref{fig:snnxgnn-mse} provides graphical representations of the ground-truth costs computed by SAA, along with their confidence intervals and the predicted cost values from the estimation methods. The light yellow area represents the cost provided by SAA wrapped by its 95\% confidence interval. The darker continuous curve corresponds to \texttt{SNN}, and the lighter continuous curve corresponds to the \texttt{GNN}. To facilitate the visualization, we display the results for a random subset of 300 districts ordered by increasing ground-truth cost. The figure displays these graphs for two different metropolitan areas (Bristol and London), considering centered depots and three possible configurations for the $n$ and $t$ factors. Generally, we can observe that \texttt{GNN} provides better estimates than \texttt{SNN} for each scenario.
Beyond this, comparing the figures with $t=3$ and $t=12$, we notice that lower $t$ values lead to much smaller errors, i.e., \texttt{GNN} and \texttt{SNN} estimates are much closer to the true cost. Larger values of $n$ also tend to diminish prediction accuracy, but this effect is less marked than when changing the district sizes.\\

\begin{figure}[htbp]
    \centering
    \includegraphics[width=0.9\linewidth]{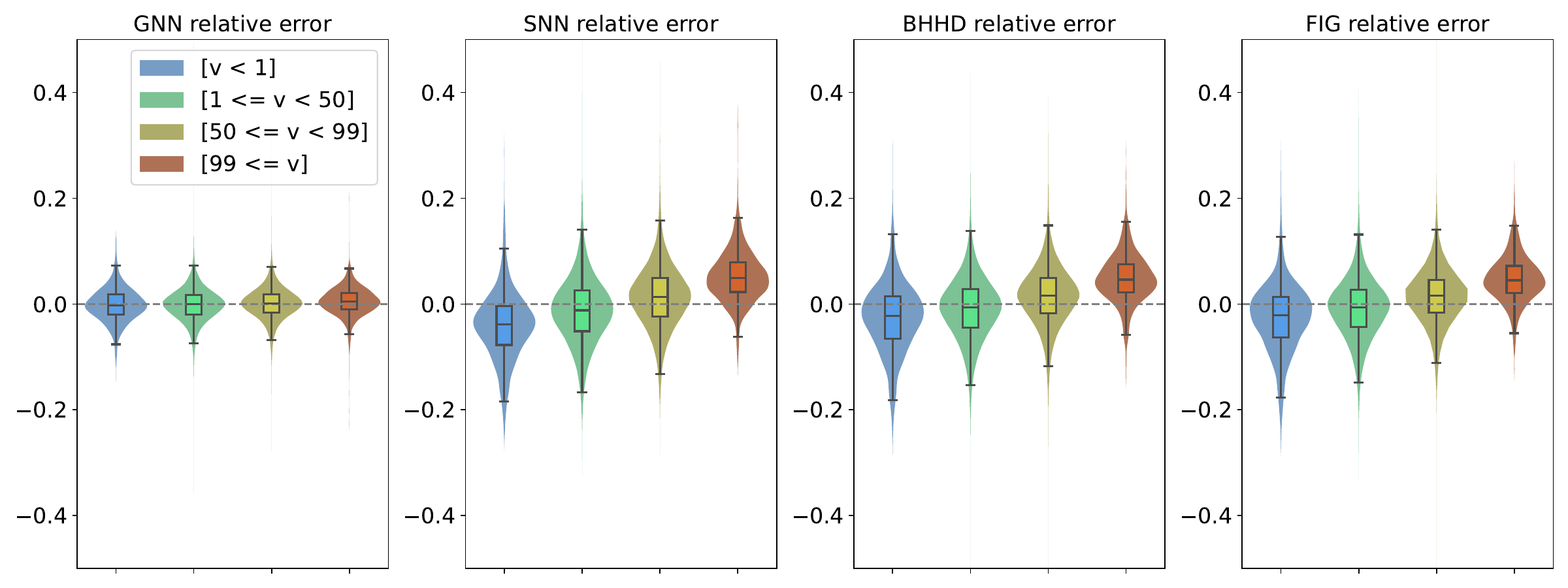}
    \caption{Performance analysis of the different estimators for districts with varying compactness}
    \label{fig:cost-all-compactness}
\end{figure}

\begin{figure}[htbp]
    \centering
    \includegraphics[width=0.9\linewidth]{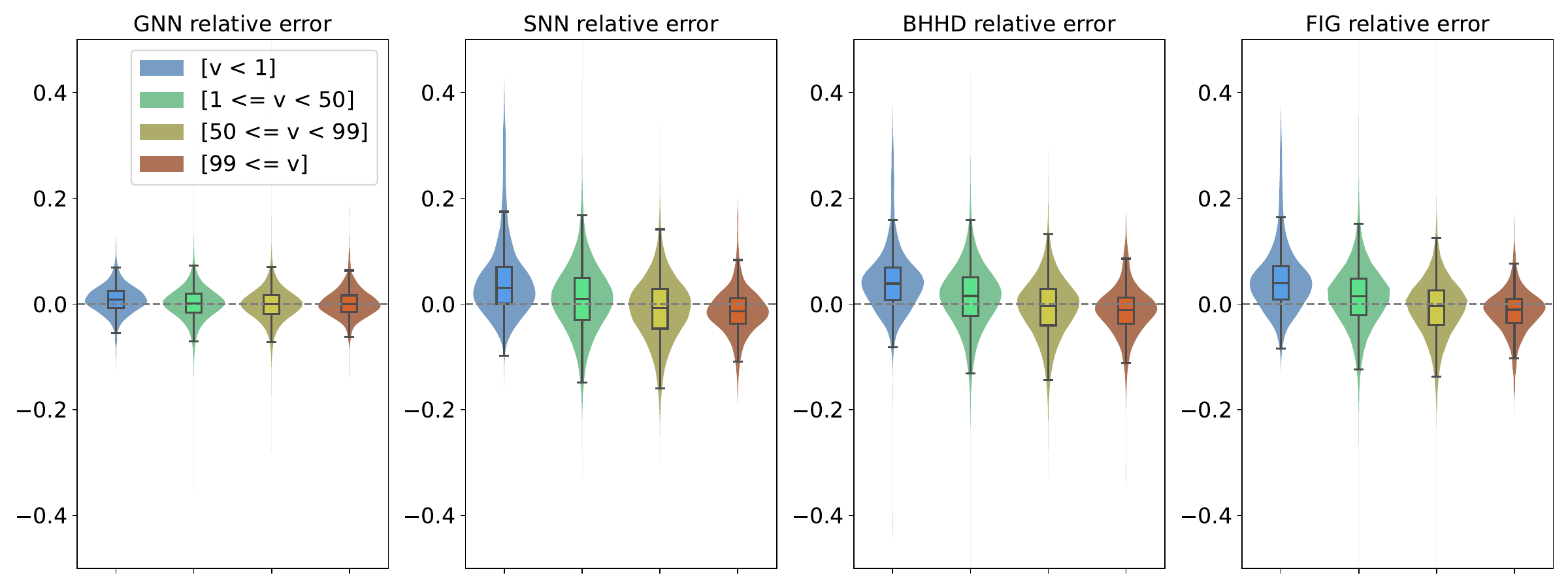}
    \caption{Performance analysis of the different estimators for districts with varying elongation}
    \label{fig:cost-all-elongation}
\end{figure}

\noindent
\textbf{Impact of district geometry on estimation performance.} 
This experiment evaluates how district geometry affects the performance of prediction models. We focus on two geometrical characteristics: compactness (defined earlier) and elongation toward the depot (defined below). Let $\sigma^{\texttt{N}}_d$ and $\sigma^{\texttt{O}}_d$ represent the standard deviations of customer positions in district $d$, measured along the axis between the depot and the district barycenter, and the perpendicular axis, respectively. The ratio $\sigma^{\texttt{O}}_d/\sigma^{\texttt{N}}_d$ indicates the district's elongation: a high ratio reflects elongation perpendicular to the depot, while a low ratio suggests elongation toward the depot. To avoid districts that are near or include the depot, we retain only the 50\% most distant districts from the test set. These districts are then split into four categories based on their compactness (or elongation) percentile value $v$ ($[v<1]$, $[1 \leq v<50]$, $[50 \leq v<99]$, $[v \geq 99]$). For each method, we report the relative cost estimation error of districts within each category, as shown in Figure~\ref{fig:cost-all-compactness} for compactness and Figure~\ref{fig:cost-all-elongation} for elongation, using violin plots.

As shown in this experiment, the simpler baselines (\texttt{BHHD}, \texttt{FIG}, \texttt{SNN}) tend to underestimate costs in less compact districts and overestimate them in more compact ones. As a result, they often favor solutions with overly compact districts, although, as noted in Section~\ref{subsec:impact-on-districting}, this may not always be optimal for districting and routing. A similar trend is observed with the elongation metric: districts elongated toward the depot (from where the delivery vehicles originate) are overestimated, while those elongated perpendicularly are underestimated. In contrast, the GNN produces consistent prediction errors across all metrics and categories.\\

\noindent
\textbf{Ablation analysis of the GNN features.}
We selected BU features commonly used in previous methods to train the GNN. Features like area and population were included both directly and as square roots, as seen in certain approximation formulas (e.g., BHHD). Although not exhaustive or minimalist, this approach leverages the neural network's ability to identify and reinforce useful features, reducing the need for extensive feature engineering. 
To provide further insight into the importance of each feature, Table~\ref{tab:ablation} summarizes the results of an ablation analysis, where the GNN was trained and tested with different subsets of the features. In this analysis, we retained the membership features defining the district for the prediction task, while selectively adding or removing one feature.

\begin{table}[htbp]
    \setlength\tabcolsep{8pt}
    \centering
        \caption{Ablation Analysis: Prediction error when removing or adding GNN features}\label{tab:ablation}
    \scalebox{0.85}
    {
        \begin{tabular}{lcc@{\hspace*{1.5cm}}lc}
        \toprule
RMSE with all features:&3.14&& RMSE with membership only:&4.67\\
\midrule
$-$ population&3.19&&$+$ population&4.06\\
$-$ sqrt-population&3.24&&$+$  sqrt-population&4.11\\
$-$ area&3.27&&$+$  area&3.73\\
$-$ sqrt-area&3.21&&$+$  sqrt-area&3.60\\
$-$ density&3.27&&$+$  density&3.74\\
$-$ perimeter&3.22&&$+$  perimeter&3.57\\
$-$ distance to the depot&3.24&&$+$  distance to the depot&3.94\\
$-$ pop \& sqrt-pop&3.28&&&\\
$-$ area \& sqrt-area&3.29&&&\\
\bottomrule
    \end{tabular}
    }
\end{table}

Since the GNN is designed for training and application in a specific city, the membership variables alone enable the model to extrapolate known district costs from the training set to unseen districts, resulting in an average RMSE of 4.67. However, incorporating additional features significantly enhances this interpolation, reducing the prediction error to an RMSE of 3.14 when all features are used. As shown in the table, the individual impact of each feature is not immediately noticeable when deactivated, as the remaining features often compensate for the missing data. Still, when each feature is added to the membership information, it provides valuable insights that the GNN can leverage. The perimeter and square root of the area contribute the most to improving accuracy, reducing RMSE from 4.67 to 3.57 and 3.60, respectively. The importance of the square root of the area is expected, as it is a key parameter in approximation formulas like BHHD.

\subsection{Analysis -- Operational Efficiency of the Districting-and-Routing Solutions}
\label{subsec:impact-on-districting}

This section evaluates the influence of various cost estimation methods within the context of a districting-and-routing solution approach. We conduct our experiments with four estimation methods: \texttt{BHHD}, \texttt{FIG}, \texttt{SNN}, and \texttt{GNN}.
The iterated local search metaheuristic (see Section~\ref{sec:districting-algorithm}) has been designed to be generic and accommodate any cost-estimation algorithm. This versatility enables us to switch estimation approaches and gauge their effects on the districting-and-routing solutions. We refer to the variations of the ILS employing different cost estimators as \texttt{GNN-ILS}, \texttt{BHHD-ILS}, \texttt{FIG-ILS}, and \texttt{SNN-ILS}. We systematically evaluate the district quality using SAA in conjunction with LKH across 500 test scenarios to obtain cost estimates. All methods are run on a single CPU thread and have the same time budget available: 180, 600, and 1200 seconds for instances of 60, 90, and 120 BUs, respectively. \myblue{We note that this moderate time budget, up to 20 minutes in the largest case, is sufficient for convergence of the ILS (with optimality gaps below $0.1\%$) for each cost oracle, as shown in our calibration and validation experiments in the Appendix. Consequently, extending the runs yields no tangible benefit when fast district-cost predictors are used. Besides this, maintaining short run times is desirable in practice, as it enables rapid ``what-if'' analyses (e.g., considering different district sizes) and iterative decision support (e.g., reviewing, adjusting, or partially fixing districts) with a human in the loop, without retraining or long computation delays.}
\\

\noindent
\textbf{Comparisons with full-knowledge solutions on smaller instances.}
For the smallest instances, it is possible to enumerate all the possible districts, evaluate their costs using the SAA approach on the test scenarios, and solve the graph partitioning formulation of Equations~(\ref{set-part-equation:objective}--\ref{set-part-equation:d-as-binary}). This approach gives us full-knowledge solutions that can be used as a bound on the best possible cost ($z_\textsc{saa-test}$). However, this requires a massive computational effort. As an illustrative example, let us consider an instance with $n=60$ and $t=6$. On average, approximately 700,000 feasible connected districts satisfy the problem's constraints. Consequently, $500 \times 700,000$ TSPs must be solved (one for each potential district and test scenario) to estimate ground-truth costs and subsequently solve the related set-partitioning formulation. Given that each TSP solution requires roughly 0.1 seconds, the total computational effort for the full-knowledge optimal approach is roughly 405 CPU-core days \emph{for each instance}. We could undertake this computation only for instances with $(n,t) \in \{(60,3),(60,6),(90,3),(90,6),(120,3)\}$ using a large computer grid. We highlight that any search method relying on SAA encounters similar limitations, as solution-cost evaluation is a time-intensive process. Relying on a supervised learning approach to estimate the costs alleviates this issue.

Table~\ref{tab: gap-ils-optimal} therefore compares the cost $z$ of the \texttt{ILS} solutions obtained with different district-cost estimation methods (using out-of-sample cost evaluations on the test set) with that of the full-knowledge solutions and reports the results as percentage gap: $\text{Gap}(\%) = 100 \times (z - z_\textsc{saa-test})/z_\textsc{saa-test}$ for each instance. The results are presented separately and averaged over the different instance sizes and depot-location configurations. As seen in these results, the proposed $\texttt{GNN-ILS}$ consistently achieves the smallest error gap relative to the full-knowledge solutions, outperforming all the methods with other estimators. The differences in performance appear to be more marked on instances with a larger number of BUs ($n$) and with larger districts ($t$) and already reach $7.48\%$ and $8.41\%$ when $t = 6$.\\

\begin{table}[!htbp]
    \vspace*{0.3cm}
    \setlength\tabcolsep{6.5pt}
    \centering
        \caption{Gap (\%) of the \texttt{ILS} with different district-cost estimation methods relative to full-knowledge solutions}\label{tab: gap-ils-optimal}
    \scalebox{0.85}
    {
        \begin{tabular}{ccccccccccc}\toprule
            \multirow{2}{*}{$\bm{n}$} & \multirow{2}{*}{$\bm{t}$} & \multicolumn{4}{c}{\textbf{Central Depot \{\texttt{C}\}}} & \multicolumn{4}{c}{\textbf{Off-Centered Depot \{\texttt{NE},\texttt{NW},\texttt{SE},\texttt{SW}\}}} \\
            \cmidrule(lr){3-6} \cmidrule(lr){7-10}
            
            & & \texttt{GNN-ILS} & \texttt{BHHD-ILS} & \texttt{FIG-ILS}  & \texttt{SNN-ILS} & \texttt{GNN-ILS} & \texttt{BHHD-ILS} & \texttt{FIG-ILS} & \texttt{SNN-ILS} \\
            \midrule
            \multirow{2}{*}{60} 
            &3 &\textbf{1.86} &6.16 &4.71 &5.28 &\textbf{1.40} &3.18 &3.14 &3.03 \\
            &6 &\textbf{0.95} &10.82 &10.92 &13.26 &\textbf{1.07} &8.39 &8.24 &7.90 \\
            
            \midrule
            
            \multirow{2}{*}{90} 
            &3 &\textbf{2.36} &5.49 &5.58 &5.39 &\textbf{1.37} &2.98 &2.99 &3.02 \\
            &6 &\textbf{2.22} &10.89 &11.18 &11.04 &\textbf{1.70} &8.20 &8.41 &7.48 \\

            \midrule
            
            120 &3 &\textbf{2.45} &5.33 &4.99 &5.41 &\textbf{1.22} &2.87 &2.95 &2.85 \\
 \midrule 
 \midrule 
   \multicolumn{2}{c}{Minimum} &\textbf{0.95} & 5.33 & 4.71 & 5.28 & \textbf{1.22} & 2.87  & 2.95 & 2.85 \\
 \multicolumn{2}{l}{Average} &\textbf{1.97} &7.74 & 7.48 & 8.08 & \textbf{1.35} & 5.12  & 5.15 & 4.86 \\
  \multicolumn{2}{l}{Maximum} &\textbf{2.45} &10.89 & 11.18 & 11.04 & \textbf{1.70} & 8.39  & 8.41 & 7.90 \\
    \bottomrule
    
    \end{tabular}
    }
\end{table}

\noindent
\textbf{General comparison.}
Considering all the instances, full-knowledge solutions are, in most cases, unavailable, and therefore we present the percentage solution-cost gaps for each method (\texttt{BHHD-ILS}, \texttt{FIG-ILS}, and \texttt{SNN-ILS}) in comparison to \texttt{GNN-ILS}, employing out-of-sample evaluation on the test set. Table \ref{tab:ils-tsp-solver} details these gaps across all instances.

As seen in this experiment, the \texttt{GNN-ILS} consistently demonstrates superior performance compared to other methods, often by a significant margin. The average gap ranges between $10.12\%$ and $13.57\%$ depending on the method and depot configuration. \myblue{Interestingly, this performance gap is significantly larger than the differences in predictive errors (MAPE measured in Section~\ref{subsec:results-predictive-methods}), suggesting that optimization tends to amplify estimator errors by disproportionately selecting districts whose costs are underestimated}. Additionally, we observe that our approach's relative performance strengthens as district size increases. This result underscores the ability of the GNN to exploit the inherent structural characteristics of the districting problem, resulting in improved estimations for larger districts.

\begin{table}[!htbp]
\setlength\tabcolsep{8pt}
 \centering
  \caption{Relative difference from \texttt{BHHD-ILS}, \texttt{FIG-ILS} and \texttt{SNN-ILS} solutions compared to \texttt{GNN-ILS}.}\label{tab:ils-tsp-solver}
 \scalebox{0.85}
{
 \begin{tabular}{cccccccc}
 \toprule
 \multirow{2}{*}{$\bm{n}$} & \multirow{2}{*}{$\bm{t}$} & \multicolumn{3}{c}{\textbf{Central Depot \{\texttt{C}\}}} & \multicolumn{3}{c}{\textbf{Off-Centered Depot \{\texttt{NE},\texttt{NW},\texttt{SE},\texttt{SW}\}}} \\
 \cmidrule(lr){3-5} \cmidrule(lr){6-8}
 & & \texttt{BHHD-ILS} & \texttt{FIG-ILS} & \texttt{SNN-ILS} 
 & \texttt{BHHD-ILS} & \texttt{FIG-ILS} & \texttt{SNN-ILS} \\

 \midrule
 \multirow{5}{*}{60} 

    &3 &4.22 &\textbf{2.80} &3.37 &1.77 &1.72 &\textbf{1.62} \\
    &6 &\textbf{9.77} &12.19 &9.87 &7.25 &\textbf{6.76} &7.10 \\
    &12 &11.57 &15.36 &\textbf{10.90} &12.40 &12.77 &\textbf{12.31} \\
    &20 &17.06 &15.57 &\textbf{13.90} &12.69 &13.73 &\textbf{12.32} \\
    &30 &18.01 &\textbf{17.60} &19.31 &\textbf{11.57} &12.39 &11.59 \\
 \midrule 
 \multirow{5}{*}{90} 
 
    &3 &3.06 &3.16 &\textbf{2.96} &\textbf{1.59} &1.60 &1.62 \\
    &6 &\textbf{8.50} &8.63 &8.79 &6.45 &\textbf{5.73} &6.56 \\
    &12 &16.57 &16.29 &\textbf{14.61} &12.70 &13.55 &\textbf{12.16} \\
    &20 &22.90 &20.67 &\textbf{17.78} &18.00 &19.15 &\textbf{17.88} \\
    &30 &21.00 &21.38 &\textbf{20.86} &16.38 &\textbf{16.19} &16.98 \\

 \midrule 
 \multirow{5}{*}{120} 
 
    &3 &2.81 &\textbf{2.48} &2.88 &1.63 &1.71 &\textbf{1.61} \\
    &6 &7.48 &\textbf{7.14} &8.40 &5.56 &5.55 &\textbf{5.28} \\
    &12 &16.92 &\textbf{14.71} &15.24 &11.18 &11.29 &\textbf{10.87} \\
    &20 &\textbf{19.36} &23.27 &24.36 &15.50 &15.66 &\textbf{14.96} \\
    &30 &\textbf{20.58} &22.26 &21.13 &\textbf{17.71} &18.14 &18.88 \\
 \midrule 
 \midrule 
   \multicolumn{2}{l}{Minimum} &2.81 & \textbf{2.48} & 2.88 & \textbf{1.59} &  1.60  & 1.61  \\
 \multicolumn{2}{l}{Average} &13.32 &13.57 &\textbf{12.96} &10.16 &10.40 &\textbf{10.12} \\
 \multicolumn{2}{l}{Maximum} &\textbf{22.90} & 23.27 & 24.36 &\textbf{18.00} &19.15 &18.88 \\
 \bottomrule
 \end{tabular}
 }
 \end{table}

\noindent
\textbf{Impact of the depot location.} 
Depots are often located at the periphery of metropolitan areas in most practical situations. Therefore, we must analyze how the geographical location of the depot generally impacts the performance of the proposed cost estimation and districting methodology. As the depot is located further away, the distance needed to drive from it to the first delivery location and return to it from the last delivery location increases. We can empirically measure the share of this back-and-forth distance in the total traveled distance, as seen in Table~\ref{tab:ils-costs-depot-vs-district}, for different subsets of the instances and relate it to the performance Gap(\%) between \texttt{GNN-ILS} and \texttt{SNN-ILS} (the second best approach). As seen in this experiment, differences in performance between the two approaches tend to increase as the share of the back-and-forth distance diminishes, i.e., when the largest part of the distance comes from trips between successive customers. This confirms the intuition that GNN predictions give an edge over simpler estimation approaches in obtaining accurate estimation in situations where costs predominantly arise from intricate routing within districts instead of simpler round trips to the depot.

\begin{table}[!htbp]
\setlength\tabcolsep{10pt}
\centering
\caption{Gap (\%) between \texttt{SNN-ILS} and \texttt{GNN-ILS} and share of back-and-forth distance (\%)}
\label{tab:ils-costs-depot-vs-district}
\scalebox{0.8}
{
 \begin{tabular}{cccccc}\toprule
 \multirow{3}{*}{$\bm{t}$} & \multicolumn{2}{c}{\textbf{Central Depot \{\texttt{C}\}}} & \multicolumn{2}{c}{\textbf{Off-Centered Depot \{\texttt{NE},\texttt{NW},\texttt{SE},\texttt{SW}\}}} \\ \cmidrule(lr){2-3} \cmidrule(lr){4-5}
 & \textbf{Gap}(\%) & \textbf{Back-forth distance(\%)} & \textbf{Gap(\%)} & \textbf{Back-forth distance(\%)} \\ 
\midrule
3  & 4.90 &22.65 &1.25 &53.94 \\
6  & 17.32 &14.29 &5.47 &42.53 \\
12 & 25.24 &8.78 &14.15 &31.02 \\
20 & 18.49 &6.12 &12.54 &24.14 \\
30 & 24.99 &4.64 &15.55 &20.27 \\
\bottomrule
\end{tabular}
}
\end{table}

\begin{figure}[!htbp]

  \hspace*{1cm}
  \begin{minipage}[t]{0.37\textwidth}
    \includegraphics[width=\linewidth]{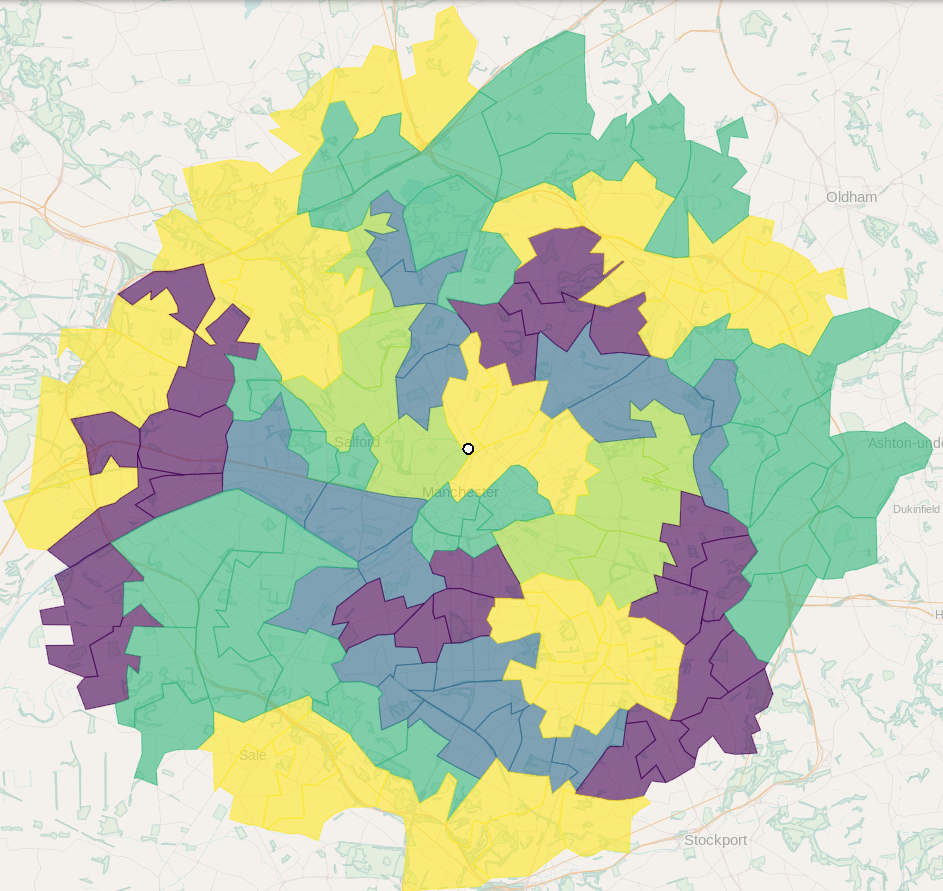}
    \begin{small} (a) \texttt{SNN-ILS} solution for \texttt{Manchester} \\ \texttt{Centered} Depot, $n=120$, $k=20$\\ Compactness Score: $0.327$,  Cost: $794.15$ \end{small}
  \end{minipage}\hspace*{2cm}
  \begin{minipage}[t]{0.37\textwidth}
     \includegraphics[width=\linewidth]{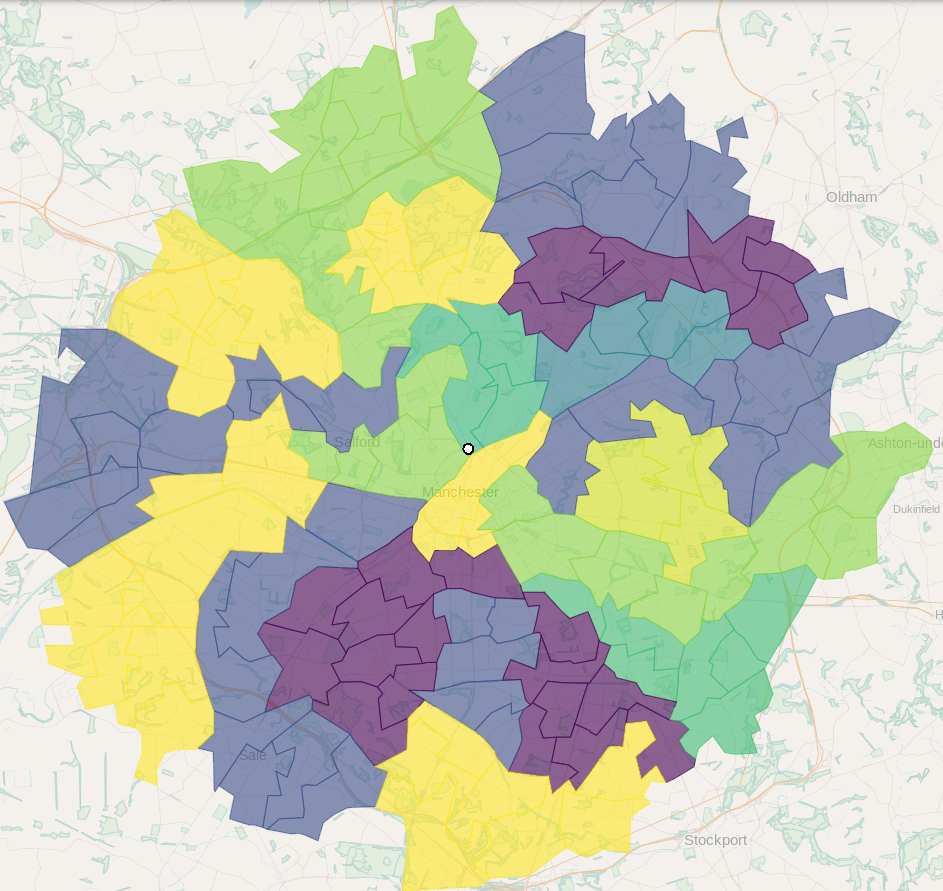}
    \begin{small} (b) \texttt{GNN-ILS} solution for \texttt{Manchester} \\ \texttt{Centered} Depot, $n=120$, $k=20$ \\ Compactness Score: $0.392$,  Cost: $714.18$ \end{small}
  \end{minipage}  
  
  \vspace{0.4cm}

  \hspace*{1cm}
  \begin{minipage}[t]{0.37\textwidth}
    \includegraphics[width=\linewidth]{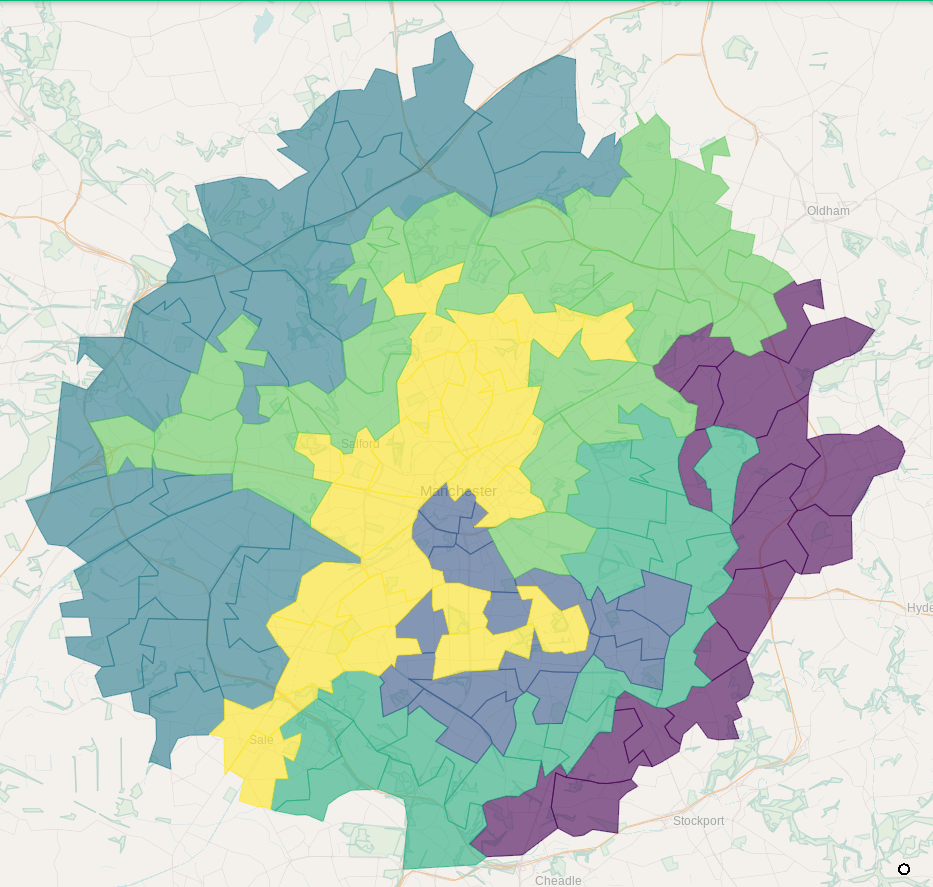}
     \begin{small} (c) \ \texttt{SNN-ILS} solution  for \texttt{Manchester} \\ \texttt{SE} Depot, $n=120$, $k=6$ \\ Compactness Score: $0.245$,  Cost: $511.19$ \end{small}
  \end{minipage}\hspace*{2cm}
  \begin{minipage}[t]{0.37\textwidth}
     \includegraphics[width=\linewidth]{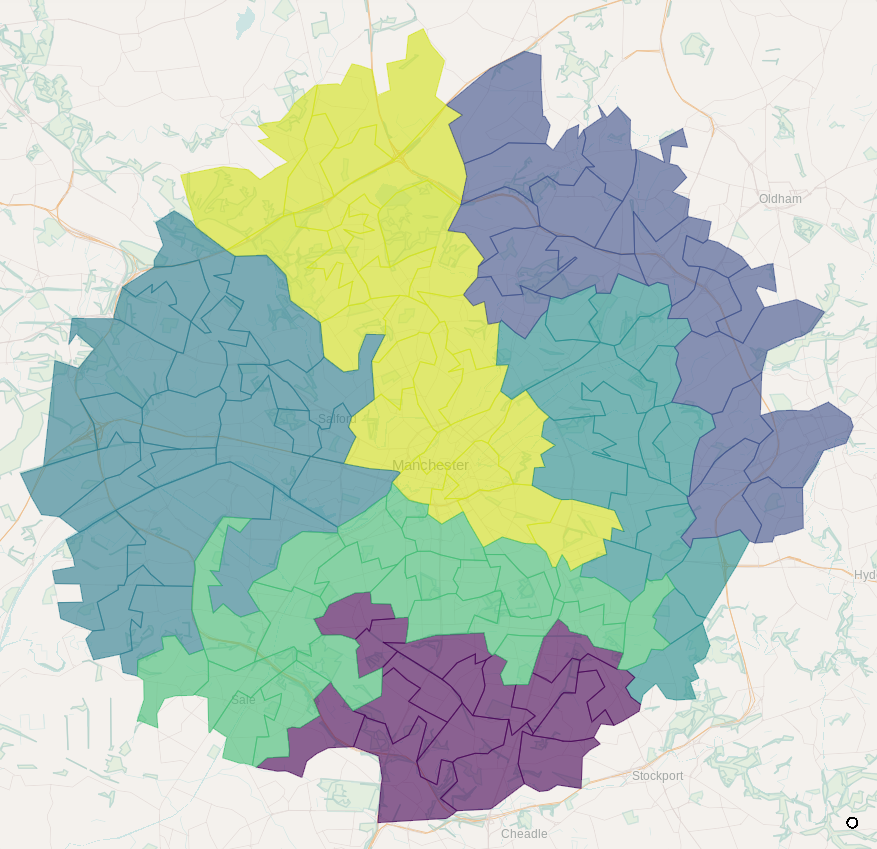}
     \begin{small} (d) \ \texttt{GNN-ILS} solution for \texttt{Manchester} \\ \texttt{SE} Depot, $n=120$, $k=6$ \\ Compactness Score: $0.360$,  Cost: $452.07$ \end{small}
  \end{minipage}
  \vspace*{0.3cm}
  
  \caption{Visual comparison of the solutions produced by \texttt{SNN-ILS} and \texttt{GNN-ILS}}
  \label{fig:ils-tsp}
\end{figure}

\noindent
\textbf{Visual analysis.}
Figure~\ref{fig:ils-tsp} compares the solutions of \texttt{GNN-ILS} and \texttt{SNN-ILS} on two instances based on the Manchester metropolitan area, the first case with $20$ districts and a centered depot, and the second case with $6$ districts and an off-centered depot in the south-eastern part of the area. Each district is represented by a distinct color, and the depot by a dot. In both cases, the operational cost of the \texttt{GNN-ILS} solution is significantly inferior to that of the \texttt{SNN-ILS} solution, with improvements of 10.08\% and 11.57\%, respectively. The solutions of the \texttt{GNN-ILS} are also (i) more compact and (ii) exhibit a geometry in which the districts are \emph{elongated towards the depot location}. As discussed in the next subsection, this latter characteristic is intuitive and 
critical to achieving good performance. None of the features of BHHD, FIG, or SNN measure this geometrical characteristic; therefore, these structures do not naturally emerge. In contrast, the GNN learns from past good combinations of BUs and has the ability to learn to promote such characteristics.

\subsection{\myblue{Discussions}}

\noindent
\textbf{Geometrical characteristics of good districts.} Compactness scores are often used as a proxy for routing efficiency within districts \citep{Benzarti2013,Bruno2021}, and Figure~\ref{fig:ils-tsp} illustrates that cost-efficient districts are generally compact. However, our results also show that operational efficiency extends beyond simple compactness measures and, therefore, is not always fully aligned with this objective.

In particular, Table \ref{tab:compactness-measure} reports the average compactness of the solutions generated using the ILS with the different estimators. When the depot is centrally located, the most efficient solutions in terms of cost (produced by the \texttt{GNN-ILS}) also feature the highest level of compactness, suggesting that compactness and routing efficiency coincide well in this case. However, this holds to a lesser extent for large instances and when the depot is off-centered, as optimal vehicle routing solutions typically exhibit routes that are elongated toward the depot rather than circular (e.g., see the solutions of the instances from \citet{Uchoa2017} accessible at \url{http://vrp.galgos.inf.puc-rio.br/index.php/en/}). \myblue{Because these elongated patterns are intrinsic to efficient routes, one can expect strategic districts to exhibit similar shapes.}

\begin{table}[!htbp]
    \setlength\tabcolsep{8pt}
    \renewcommand{\arraystretch}{1.1}
    \centering
    \caption{Average compactness of the districts produced by \texttt{GNN-ILS}, \texttt{BHHD-ILS}, \texttt{FIG-ILS}, and \texttt{SNN-ILS}}\label{tab:compactness-measure}
    \scalebox{0.85} 
    {
    \begin{tabular}{cl@{\hspace*{1cm}}cccc}
    \toprule
    \multirow{2}{*}{\textbf{Depot}} & \multirow{2}{*}{$\bm{t}$} & \multicolumn{4}{c}{\textbf{Average Compactness}} \\
    & & \texttt{GNN-ILS} & \texttt{BHHD-ILS} & \texttt{FIG-ILS} & \texttt{SNN-ILS} \\

    \midrule
    \multirow{5}{*}{\shortstack{Central Depot \\ $\{\texttt{C}$\}}} 

    &3 &\textbf{0.387} &0.375 &0.382 &0.374 \\
    &6 &\textbf{0.386} &0.341 &0.336 &0.339 \\
    &12 &\textbf{0.400} &0.308 &0.316 &0.307 \\
    &20 &\textbf{0.409} &0.317 &0.315 &0.326 \\
    &30 &\textbf{0.421} &0.350 &0.352 &0.354 \\

    \midrule
    \multirow{5}{*}{\shortstack{Off-Centered Depot \\ $\{\texttt{NE},\texttt{NW},\texttt{SE},\texttt{SW}\}$}}

    &3 &0.367 &\textbf{0.369} &\textbf{0.369} &0.368 \\
    &6 &0.322 &0.340 &\textbf{0.344} &0.338 \\
    &12 &0.293 &0.314 &\textbf{0.315} &0.313 \\
    &20 &\textbf{0.337} &0.328 &0.320 &0.326 \\
    &30 &\textbf{0.375} &0.352 &0.357 &0.355 \\

    \bottomrule
    \end{tabular}
    }
    \end{table}

\myblue{
These characteristics are clearly visible in the solutions obtained with the GNN (e.g., Figure~\ref{fig:ils-tsp}d). Interestingly, the GNN can guide the optimization towards such geometries, finding an efficient middle ground between compactness and elongation towards the depot, even without explicitly defined features measuring these properties. This underscores the flexibility of the supervised learning approach: as it gauges solution quality solely by imitation of known training districts.
Consequently, the same methodology can be readily applied in various contexts, such as different routing subproblems, additional operational constraints, or network-based distance metrics, by training on district examples labeled with a different cost without any other architectural modification.}\\

\myblue{
\noindent
\textbf{Experiments considering network-based distances.} To further illustrate the flexibility of our approach, we evaluated its performance when distances are computed on the road network, therefore accounting for other real-world complexities such as tunnels, bottlenecks, and geographical barriers. Only minimal adaptation was required for this setting: the sole modification was to use road-network distances in the TSPs when generating the training districts. For this purpose, all-pairs shortest paths between delivery locations were computed using the contraction hierarchies algorithm of \citet{Geisberger2012} and \citet{Dibbelt2016}, implemented in the \texttt{RoutingKit} library (\url{https://github.com/RoutingKit}). All the district-cost estimation models (GNN, BHHD, GIF, and SNN) were subsequently retrained or recalibrated on these new data to evaluate both their predictive accuracy and their impact on districting-and-routing performance. For brevity, these experiments were conducted on the subset of instances with $n=120$ basic units and a centrally located depot. The corresponding results are summarized in Table~\ref{table:real-distances}.}

\begin{table}[ht]
\setlength\tabcolsep{8pt}
\renewcommand{\arraystretch}{1.15}
\centering
\caption{Accuracy and operational performance of district-cost estimators under road-network distances}
\label{table:real-distances}
\scalebox{0.82}{
\begin{tabular}{cccccccccccccc}
\toprule
\multirow{2}{*}{$\bm{n}$}  & \multirow{2}{*}{$\bm{t}$} &
\multirow{2}{*}{$\phi_{\mathrm{SAA}}^{\text{test}}$} &
\multicolumn{4}{c}{\textbf{RMSE}} &
\multicolumn{4}{c}{\textbf{MAPE}} &
\multicolumn{3}{c}{\textbf{Gap(\%) vs \texttt{GNN-ILS}}} \\
\cmidrule(lr){4-7} \cmidrule(lr){8-11} \cmidrule(lr){12-14}
 &  & &
 \texttt{GNN} & \texttt{BHHD} & \texttt{FIG} & \texttt{SNN} &
 \texttt{GNN} & \texttt{BHHD} & \texttt{FIG} & \texttt{SNN} &
 \texttt{BHHD} & \texttt{FIG} & \texttt{SNN} \\
\midrule
\multirow{5}{*}{120}
 & 3  & 98.08  & \bfseries 6.03  & 19.48 & 16.08 & 15.68 & \bfseries 4.88 & 15.15 & 13.35 & 12.90 & 2.76 & 2.52 & 2.36 \\
 & 6  & 101.87 & \bfseries 5.19  & 18.62 & 14.17 & 13.86 & \bfseries 3.59 & 14.38 & 11.50 & 11.06 & 6.17 & 5.80 & 5.54 \\
 & 12 & 145.25 & \bfseries 20.12 & 43.20 & 37.28 & 36.98 & \bfseries 6.77 & 24.92 & 22.33 & 21.91 & 9.62 & 9.02 & 8.66 \\
 & 20 & 158.26 & \bfseries 18.87 & 37.31 & 30.96 & 31.27 & \bfseries 6.38 & 20.08 & 16.57 & 16.56 & 8.26 & 10.34 & 6.74 \\
 & 30 & 172.06 & \bfseries 12.65 & 30.96 & 25.26 & 25.37 & \bfseries 4.90 & 16.01 & 12.71 & 12.91 & 14.81 & 18.22 & 10.75 \\
\midrule
\multicolumn{2}{c}{\textbf{Average}} & 135.11 & \bfseries 12.57 & 29.92 & 24.75 & 24.63 & \bfseries 5.30 & 18.11 & 15.29 & 15.07 & 8.32 & 9.18 & 6.81 \\
\bottomrule
\end{tabular}
}
\end{table}

\myblue{
The results on road-network data follow the same overall trends observed with Euclidean distances. The GNN achieves the lowest RMSE and MAPE on hold-out district-cost test sets, followed by the SNN, FIG, and BHHD estimators. In terms of operational efficiency, the districts produced by \texttt{GNN-ILS} remain, on average, about 7\% to 9\% more cost-effective than those obtained with alternative baselines.\\
}

\noindent
\textbf{\myblue{Other comparisons with reference solutions.}} \myblue{As mentioned earlier, using a full SAA approach that solves multiple TSPs to evaluate district costs at each step would be prohibitively slow. As an alternative, we implemented another simplified baseline that considers a single demand realization scenario. With this method, each district cost evaluation reduces to solving a deterministic TSP.
In our experiments, this approach underperformed the GNN-ILS by 1.88\% and 1.82\% on average for Euclidean instances with central and off-center depots, respectively, and by 2.71\% in the road-network cases, highlighting the value of leveraging stochastic information. Nevertheless, it performed favorably compared to ILS with the BHHD, FIG, or SNN estimators, likely because, unlike CA formulas relying on district statistics, the single-scenario TSP approach does not systematically steer the districting method towards certain suboptimal geometrical characteristics such as compactness. Its main drawback, however, is computational: even a single TSP solution per move evaluation in a districting method is considerably slower than evaluating a CA formula or performing neural network inference, permitting very few local search descents for districting optimization over a given time.

We also investigated replacing the heuristic TSP evaluations based on the Lin–Kernighan algorithm with machine-learning-based estimators such as the approach proposed by \citet{Varol2024}, which predicts tour lengths from given point coordinates. However, under our experimental configuration with tours of around 100 customers, the fast parametrization of LKH consistently finds near-optimal solutions (optimal in 92\% of the cases and within 0.1\% of the exact solution on average) within 0.1 seconds, whereas the point-based estimator requires complex feature generation and typically yields larger TSP-estimation errors (around 1\%) at a slower runtime of 0.6 to 1.2 seconds. Consequently, current point-based TSP estimators provide no computational or accuracy advantage over LKH. By contrast, the trained GNN takes roughly 8ms for inference on a given district, and takes demand distribution information as input instead of being tied to the inherent variability of a single demand scenario.}

\myblue{Finally, to broaden our comparison, we evaluated a constrained version of K-means (CK-means) using the algorithm of \citet{Levy-Kramer_k-means-constrained_2018}. This method is extremely fast, requiring only a few seconds and no parameter calibration, but its solutions were, on average, 4.46\% and 3.31\% worse than those obtained with \texttt{GNN-ILS} in the Euclidean and road-network cases, respectively. Moreover, since the algorithm does not guarantee spatial connectivity, it frequently produced infeasible district configurations}.\\

\noindent
\textbf{Generalization capabilities.}
Building on the architecture proposed in this work, a subsequent study by \citet{Ahmed2024} also extended the methodology into an end-to-end, decision-aware learning framework that integrates a combinatorial optimization layer within the GNN pipeline. In this formulation, the network parametrizes a surrogate problem based on a capacitated minimum spanning tree and is trained using a Fenchel-Young loss with perturbation-based gradient estimation following \citet{Dalle2022}, allowing the model to optimize the downstream decision objective directly rather than minimizing district-cost prediction error. This end-to-end approach demonstrated strong generalization across cities, across varying demographic and geographical characteristics of BUs (e.g., population density, area, and compactness), and across problem scales. It also achieved improved computational efficiency by learning on small problem instances and extrapolating to much larger ones: a model trained on instances with about 30 BUs successfully producing high-quality districting solutions for metropolitan regions containing up to 2,000 units. Overall, while the architecture introduced in this paper constitutes a fundamental building block for efficient districting-and-routing methods, subsequent advances in learning have further extended its generalization capabilities.\\

\noindent
\textbf{Summary.}
Table~\ref{tab:comparison-final} summarizes the performance, speed, and generalization properties of all the evaluated methods. The purpose of this study is not to present \texttt{ILS-GNN} as a universal tool to every districting-and-routing setting, but to highlight its robustness and the limitations of using CA formulas within optimization frameworks that tend to favor districts whose costs are most underestimated. The proposed approach provides a strong alternative whenever some preprocessing is acceptable and optimization must be repeated over multiple problem variants or partial solutions.

\begin{table}[ht]
\setlength\tabcolsep{7pt}
\renewcommand{\arraystretch}{1.3}
\centering
\caption{Summary of the considered district-cost estimation and optimization methods}\label{tab:comparison-final}
\scalebox{0.8}{
\begin{tabular}{lccHcp{0.6\linewidth}}
\toprule
\textbf{Method} & \textbf{Learning} & \textbf{Optimization} & \textbf{RMSE} & \textbf{Gap(\%)$^*$} & \textbf{Generalization to other cities} \\\midrule
ILS-BHHD & $<$1s & 20min & 4.46 & 10.79  & Retraining is optional  -- Fine-tuning regression weights for better performance requires district-cost examples  \\
ILS-FIG & $<$1s & 20min & 4.40 & 11.03  & Retraining is optional -- Fine-tuning regression weights for better performance requires district-cost examples  \\
ILS-SNN & 1s & 20min & 4.24 & 10.69  & Retraining using district-cost examples  \\
ILS-GNN & 2h & 20min &  2.31 & \textemdash{} & Retraining using district-cost examples, possible to generalize to different cities using end-to-end learning \citep{Ahmed2024} \\
ILS-TSP$^\dagger$ & \textemdash{} & 20min$^\dagger$ & \textemdash{}  & 1.83  & No learning needed \\
CK-means$^\ddagger$ & \textemdash{} & a few seconds & \textemdash{} &  4.46 & No learning needed \\
\bottomrule
\end{tabular}
}

\vspace{1mm}
\parbox{\linewidth}{
\footnotesize
\vspace*{0.05cm}
\quad $^*$ Gap (\%) relative to ILS-GNN. Euclidean Case. All districting solutions were evaluated on hold-out test scenarios. \vspace*{-0.1cm}

\quad $^\dagger$~Single-scenario evaluation using LKH3 for the current fastest and most accurate estimations of TSP costs. \vspace*{-0.1cm}

\quad $^\ddagger$~Solutions might not respect district connectivity constraints.
}
\end{table}

\section{Conclusions}
\label{sec:conclusions}

This paper has investigated a novel approach for strategic districting and routing, relying on a GNN for district cost estimation within an ILS, aiming to find a good partition of a metropolitan area into districts. The GNN is trained through a supervised learning process, using pre-computed TSP costs for various districting configurations derived from the LKH heuristic. Our experimental analyses demonstrate that this approach generates efficient districts, leading to substantial economic gains of 10.12\% on average (and exceeding 20\% in some cases) compared to other commonly used methods for routing cost estimations. They also show that compactness is an \myblue{insufficient proxy for efficiency} and that other geometric characteristics (e.g., elongated district shapes towards the depot) play a key role in determining long-term operational costs.

This research opens many promising research avenues.
A natural direction for future work is to analyze the applicability of GNN-based cost estimations to other stochastic problems that require time-intensive second-stage evaluations. One main advantage of the proposed GNN approach is its flexibility. It is driven by data (geographical features of the BUs and training examples) and does not require distributional information over customer requests. Consequently, it is fairly generic, and adapting it to other district-cost definitions (e.g., considering a different subproblem than a TSP, introducing additional constraints or network characteristics) boils down to replacing the approach used to generate labeled districts.

\myblue{
Another line of research involves further improving the learning methodology and its generalization capabilities. While the present work trains the GNN to minimize prediction error independently of the downstream optimization process, subsequent developments have shown that end-to-end, decision-aware training (using a surrogate optimization layer and a Fenchel–Young loss) can further enhance generalization across cities with different population densities, spatial structures, and district sizes~\citep{Ahmed2024}.}

\myblue{
Finally, a complementary perspective is to explore solution methods at the districting level that explicitly exploit finer properties of district costs. In many practical cases, adding an additional basic unit (BU) to a district leads to diminishing marginal increases in cost, suggesting an approximately submodular behavior of the district-cost function that could be leveraged algorithmically. However, submodularity violations might also arise when the inclusion of new BUs induces threshold effects, such as capacity or time-window constraints in the underlying VRP subproblems that require an additional vehicle, multi-depot settings where the optimal depot assignment changes, or spatial discontinuities introduced by bridges, tunnels, or other natural barriers. Investigating different optimization paradigms, approximation algorithms, or local-search heuristics that take advantage of specific structures to filter moves could provide valuable insights to improve solution methods, complementing the learning-based approaches examined in this study.}

\ACKNOWLEDGMENT{This project has been partly enabled by computational infrastructure provided by Calcul Québec and Compute Canada. This support is gratefully acknowledged.}

\begin{APPENDICES}

\vspace*{0.8cm}

\section{Calibration and Validation of the Iterated Local Search}
\label{appendix:ils-tuning}

To calibrate our ILS, we relied on instances that could be solved to optimality using the set partitioning formulation, therefore giving us a solid baseline. We first calibrated the perturbation parameter of the ILS-GNN, considering values of $p_{\textsc{rm}}$ among $0.500$, $0.250$, $0.150$, $0.100$, $0.050$, $0.025$, $0.020$, $0.015$, $0.010$, and $0.005$. We observed that allowing $180$ seconds of CPU time for the instances with $60$ BUs, $600$ seconds for the instances with $90$ BUs, and $1200$ seconds for the instances with $120$ BUs was sufficient to achieve near-optimal results in all cases.
Table \ref{tab: p-tuning-complete} reports the error gap of ILS-GNN for the different values of $p_{\textsc{rm}}$, measured as percentage deviation between the heuristic solutions produced by the ILS and the optimal solutions found by the set partitioning formulation on each subset of the instances. According to this experiment, we selected the best method configuration with $p_{\textsc{rm}}= 0.015$. It is important to note that with this parameter setting, the error gap of the heuristic is so small that using it instead of an exact method on the considered problem instances leads to differences that are almost imperceptible.

\begin{table}[!htbp]
\setlength\tabcolsep{5pt}
\centering
\caption{Gap (\%) of ILS-GNN to known optimal solutions for this cost oracle, for different values of~$p_{\textsc{rm}}$}
\label{tab: p-tuning-complete}
\scalebox{0.90}
{
 \begin{tabular}{lrrrrrrrrrrrr}\toprule
 & & \multicolumn{10}{c}{$\bm{p}_{\textsc{rm}}$} \\
 $\bm{n}$ & $\bm{t}$ & \textbf{0.500} & \textbf{0.250} & \textbf{0.150} & \textbf{0.100} & \textbf{0.050} & \textbf{0.025} & \textbf{0.020} & \textbf{0.015} & \textbf{0.010} & \textbf{0.005} \\
 \midrule
 \multirow{2}{*}{60}
 
 &3	&0.277	&0.186	&0.148	&0.086	&0.081	&0.087	&0.053	&0.045	&0.052	&0.049\\ 
 &6	&0.516	&0.377	&0.324	&0.224	&0.204	&0.183	&0.206	&0.184	&0.171	&0.183\\ 
 
 \multirow{2}{*}{90} 
 
 &3	&0.460	&0.337	&0.247	&0.212	&0.173	&0.155	&0.143	&0.126	&0.132	&0.136\\ 
 &6	&0.995	&0.711	&0.600	&0.480	&0.422	&0.396	&0.349	&0.317	&0.396	&0.390\\ 
 
 \multirow{2}{*}{120} 
 
 &3	&0.486	&0.386	&0.272	&0.261	&0.229	&0.199	&0.180	&0.182	&0.183	&0.192\\ 
 &6	&1.173	&0.847	&0.719	&0.610	&0.510	&0.415	&0.467	&0.514	&0.465	&0.430\\ 
 
 \midrule
 \multicolumn{2}{c}{Average} &0.651	&0.474	&0.385	&0.312	&0.270	&0.239	&0.233	&\textbf{0.228}	&0.233	&0.230\\
 \bottomrule
\end{tabular}
}
\end{table}

Moreover, since an important goal of our experiments is to analyze the characteristics of near-optimal solutions produced using the different cost-estimation oracles, we also wanted to ensure that our use of an ILS-based districting heuristic does not introduce biases that would favor certain cost oracles, which would not exist otherwise if an exact districting method had been used. Therefore, we conducted a similar analysis as for our calibration, this time measuring the error gap of the heuristic relative to the exact method for each of the different cost-estimation oracles and keeping $p_{\textsc{rm}} = 0.015$. The results of this additional analysis are reported in Table~\ref{tab:gap-comparison}, along with the total number of iterations performed by the different methods. From this analysis, we observe that the error gap of all methods remains below $0.228 \%$ on average. Moreover, the use of a heuristic does not favor the GNN-based oracle. In fact, a small opposite effect happens since the use of the GNN oracle introduces a slight overhead in the move evaluations in comparison to closed-form approximation formulas, leading to fewer iterations before termination (number of perturbations and LS loops in Line~4 of Algorithm~1 in the main paper) and consequently a slightly larger error gap. Nevertheless, it is important to note that this difference in the performance of the heuristic does not hinder the ability of the ILS-GNN to lead to districting solutions that are far superior to the others in terms of true districting costs as evaluated on the test scenarios (see Table~7 in the main paper).

\begin{table}[!htbp]
    \setlength\tabcolsep{5pt}
    \centering
        \caption{Gap (\%) and number of iterations for each cost oracle}
    \label{tab:gap-comparison}
    \scalebox{0.92}
    {
    \begin{tabular}{lrrrrrrrrrr}\toprule
    & & \multicolumn{2}{c}{\textbf{ILS-GNN}} & \multicolumn{2}{c}{\textbf{ILS-FIG}} & \multicolumn{2}{c}{\textbf{ILS-SNN}} & \multicolumn{2}{c}{\textbf{ILS-BHHD}}\\ 
    \cmidrule(lr){3-4} \cmidrule(lr){5-6} \cmidrule(lr){7-8} \cmidrule(lr){9-10}
    $\bm{n}$ & $\bm{t}$ & \textbf{Gap} & \textbf{Iter.} & \textbf{Gap} & \textbf{Iter.} & \textbf{Gap} & \textbf{Iter.} & \textbf{Gap} & \textbf{Iter.} \\
    
    \midrule
    
    \multirow{2}{*}{60} 
    
    &3 &0.045 &55920 &0.006 &319840 &0.007 &505234 &0.006 &347105 \\
    &6 &0.184 &4813 &0.011 &50962 &0.014 &76507 &0.009 &51080 \\
    
    \multirow{2}{*}{90} 
    
    &3 &0.126 &47212 &0.030 &334571 &0.029 &552337 &0.025 &364744 \\
    &6 &0.317 &7023 &0.062 &62347 &0.058 &104712 &0.063 &72002 \\
    
    \multirow{2}{*}{120} 
    
    &3 &0.182 &45447 &0.045 &414579 &0.044 &718792 &0.044 &469031 \\
    &6 &0.514 &7645 &0.091 &107971 &0.106 &180756 &0.080 &120631 \\
    
    \midrule
    
    \multicolumn{2}{c}{Average} &0.228 &28010 &0.041 &215045 &0.043 &356390 &0.038 &237432 \\
    \bottomrule
   
    \end{tabular}
    
    }
    \end{table}
 
 \textbf{}

\end{APPENDICES}


\begin{thebibliography}{61}
\providecommand{\natexlab}[1]{#1}
\providecommand{\url}[1]{\texttt{#1}}
\providecommand{\urlprefix}{URL }

\bibitem[{Ahmed et~al.(2024)Ahmed, Forel, Parmentier, \protect\BIBand{} Vidal}]{Ahmed2024}
Ahmed C, Forel A, Parmentier A, Vidal T, 2024 \emph{Districtnet: Decision-aware learning for geographical districting}. \emph{Proceedings of the 38th International Conference on Neural Information Processing Systems}, NeurIPS'24.

\bibitem[{Akkerman \protect\BIBand{} Mes(2025)}]{Akkerman2022}
Akkerman F, Mes M, 2025 \emph{{Distance approximation to support customer selection in vehicle routing problems}}. \emph{Annals of Operations Research} 350(1):269--297.

\bibitem[{Ansari et~al.(2018)Ansari, Basdere, Li, Ouyang, \protect\BIBand{} Smilowitz}]{Ansari2018b}
Ansari S, Basdere M, Li X, Ouyang Y, Smilowitz K, 2018 \emph{{Advancements in continuous approximation models for logistics and transportation systems: 1996-2016}}. \emph{Transportation Research Part B: Methodological} 107:229--252.

\bibitem[{Baddeley(2006)}]{Baddeley2006}
Baddeley A, 2006 \emph{{Spatial point processes and their applications}}. {W Weil}, ed., \emph{Stochastic Geometry}, volume 1892 of \emph{Lecture Notes in Mathematics}, 1--75 (Springer).

\bibitem[{Beardwood, Halton, \protect\BIBand{} Hammersley(1959)}]{Beardwood1959}
Beardwood J, Halton JH, Hammersley JM, 1959 \emph{{The shortest path through many points}}. \emph{Mathematical Proceedings of the Cambridge Philosophical Society} 55(9):299--327.

\bibitem[{Benzarti, Sahin, \protect\BIBand{} Dallery(2013)}]{Benzarti2013}
Benzarti E, Sahin E, Dallery Y, 2013 \emph{{Operations management applied to home care services: Analysis of the districting problem}}. \emph{Decision Support Systems} 55(2):587--598.

\bibitem[{Bozkaya, Erkut, \protect\BIBand{} Laporte(2003)}]{Bozkaya2003}
Bozkaya B, Erkut E, Laporte G, 2003 \emph{{A tabu search heuristic and adaptive memory procedure for political districting}}. \emph{European Journal of Operational Research} 144(1):12--26.

\bibitem[{Bruno et~al.(2021)Bruno, Cavola, Diglio, Laporte, \protect\BIBand{} Piccolo}]{Bruno2021}
Bruno G, Cavola M, Diglio A, Laporte G, Piccolo C, 2021 \emph{{Reorganizing postal collection operations in urban areas as a result of declining mail volumes–A case study in Bologna}}. \emph{Journal of the Operational Research Society} 72(7):1591--1606.

\bibitem[{Cappart et~al.(2023)Cappart, Chételat, Khalil, Lodi, Morris, \protect\BIBand{} Velickovic}]{cappart2021combinatorial}
Cappart Q, Chételat D, Khalil EB, Lodi A, Morris C, Velickovic P, 2023 \emph{Combinatorial optimization and reasoning with graph neural networks}. \emph{Journal of Machine Learning Research} 24(130):1--61.

\bibitem[{{\c{C}}avdar \protect\BIBand{} Sokol(2015)}]{Cavdar2015b}
{\c{C}}avdar B, Sokol J, 2015 \emph{{A distribution-free TSP tour length estimation model for random graphs}}. \emph{European Journal of Operational Research} 243(2):588--598.

\bibitem[{Chien(1992)}]{Chien1992}
Chien TW, 1992 \emph{{Operational estimators for the length of a traveling salesman tour}}. \emph{Computers and Operations Research} 19(6):469--478.

\bibitem[{Daganzo(1984)}]{daganzo1984distance}
Daganzo CF, 1984 \emph{The distance traveled to visit {N} points with a maximum of {C} stops per vehicle: An analytic model and an application}. \emph{Transportation Science} 18(4):331--350.

\bibitem[{Dai, Dai, \protect\BIBand{} Song(2016)}]{dai2016discriminative}
Dai H, Dai B, Song L, 2016 \emph{Discriminative embeddings of latent variable models for structured data}. \emph{Proceedings of The 33rd International Conference on Machine Learning}, 2702--2711.

\bibitem[{Dai et~al.(2017)Dai, Khalil, Zhang, Dilkina, \protect\BIBand{} Song}]{dai2017learning}
Dai H, Khalil EB, Zhang Y, Dilkina B, Song L, 2017 \emph{Learning combinatorial optimization algorithms over graphs}. \emph{Proceedings of the 31st International Conference on Neural Information Processing Systems}, 6351–6361.

\bibitem[{Dalle et~al.(2022)Dalle, Baty, Bouvier, \protect\BIBand{} Parmentier}]{Dalle2022}
Dalle G, Baty L, Bouvier L, Parmentier A, 2022 \emph{{Learning with combinatorial optimization layers: A probabilistic approach}}. Technical report, ArXiv: 2207.13513, \urlprefix\url{http://arxiv.org/abs/2207.13513}.

\bibitem[{Dibbelt, Strasser, \protect\BIBand{} Wagner(2016)}]{Dibbelt2016}
Dibbelt J, Strasser B, Wagner D, 2016 \emph{{Customizable contraction hierarchies}}. \emph{ACM Journal of Experimental Algorithmics} 21:1--49.

\bibitem[{Drakuli\'{c}, Michel, \protect\BIBand{} Andreoli(2025)}]{drakulic2025goal}
Drakuli\'{c} D, Michel S, Andreoli JM, 2025 \emph{{GOAL: A generalist combinatorial optimization agent learner}}. \emph{International Conference on Learning Representations}, 52465--52488.

\bibitem[{Drexl \protect\BIBand{} Schneider(2015)}]{Drexl2015}
Drexl M, Schneider M, 2015 \emph{{A survey of variants and extensions of the location-routing problem}}. \emph{European Journal of Operational Research} 241(2):283--308.

\bibitem[{Dyer \protect\BIBand{} Frieze(1985)}]{Dyer1985}
Dyer M, Frieze A, 1985 \emph{On the complexity of partitioning graphs into connected subgraphs}. \emph{Discrete Applied Mathematics} 10(2):139--153.

\bibitem[{Fey \protect\BIBand{} Lenssen(2019)}]{fey2019fast}
Fey M, Lenssen JE, 2019 \emph{Fast graph representation learning with {PyTorch Geometric}}. \emph{ICLR Workshop on Representation Learning on Graphs and Manifolds}.

\bibitem[{Figliozzi(2007)}]{Figliozzi2007}
Figliozzi MA, 2007 \emph{{Analysis of the efficiency of urban commercial vehicle tours: Data collection, methodology, and policy implications}}. \emph{Transportation Research Part B: Methodological} 41(9):1014--1032.

\bibitem[{Franceschetti, Jabali, \protect\BIBand{} Laporte(2017)}]{Franceschetti2017}
Franceschetti A, Jabali O, Laporte G, 2017 \emph{{Continuous approximation models in freight distribution management}}. \emph{TOP} 25(3):413--433.

\bibitem[{Galv{\~{a}}o et~al.(2006)Galv{\~{a}}o, Novaes, {Souza De Cursi}, \protect\BIBand{} Souza}]{Galvao2006}
Galv{\~{a}}o LC, Novaes AG, {Souza De Cursi} JE, Souza JC, 2006 \emph{{A multiplicatively-weighted Voronoi diagram approach to logistics districting}}. \emph{Computers and Operations Research} 33(1):93--114.

\bibitem[{Garc{\'{i}}a-Ayala et~al.(2016)Garc{\'{i}}a-Ayala, Gonz{\'{a}}lez-Velarde, R{\'{i}}os-Mercado, \protect\BIBand{} Fern{\'{a}}ndez}]{Garcia-Ayala2016}
Garc{\'{i}}a-Ayala G, Gonz{\'{a}}lez-Velarde JL, R{\'{i}}os-Mercado RZ, Fern{\'{a}}ndez E, 2016 \emph{{A novel model for arc territory design: Promoting Eulerian districts}}. \emph{International Transactions in Operational Research} 23(3):433--458.

\bibitem[{Geisberger et~al.(2012)Geisberger, Sanders, Schultes, \protect\BIBand{} Vetter}]{Geisberger2012}
Geisberger R, Sanders P, Schultes D, Vetter C, 2012 \emph{{Exact routing in large road networks using contraction hierarchies}}. \emph{Transportation Science} 46(3):388--404.

\bibitem[{Glorot, Bordes, \protect\BIBand{} Bengio(2011)}]{glorot2011deep}
Glorot X, Bordes A, Bengio Y, 2011 \emph{Deep sparse rectifier neural networks}. \emph{Proceedings of the Fourteenth International Conference on Artificial Intelligence and Statistics}, 315--323.

\bibitem[{Hamilton, Ying, \protect\BIBand{} Leskovec(2017)}]{hamilton2017inductive}
Hamilton WL, Ying R, Leskovec J, 2017 \emph{Inductive representation learning on large graphs}. \emph{Proceedings of the 31st International Conference on Neural Information Processing Systems}, 1025–1035.

\bibitem[{Helsgaun(2000)}]{Helsgaun2000}
Helsgaun K, 2000 \emph{{An effective implementation of the Lin-Kernighan traveling salesman heuristic}}. \emph{European Journal of Operational Research} 126(1):106--130.

\bibitem[{Holland et~al.(2017)Holland, Levis, Nuggehalli, Santilli, \protect\BIBand{} Winters}]{Holland2017}
Holland C, Levis J, Nuggehalli R, Santilli B, Winters J, 2017 \emph{{UPS optimizes delivery routes}}. \emph{Interfaces} 47(1):8--23.

\bibitem[{Horn, Hampton, \protect\BIBand{} Vandenberg(1993)}]{Horn1993}
Horn DL, Hampton CR, Vandenberg AJ, 1993 \emph{{Practical application of district compactness}}. \emph{Political Geography} 12(2):103--120.

\bibitem[{Joshi et~al.(2022)Joshi, Cappart, Rousseau, \protect\BIBand{} Laurent}]{joshi2022learning}
Joshi CK, Cappart Q, Rousseau LM, Laurent T, 2022 \emph{Learning the travelling salesperson problem requires rethinking generalization}. \emph{Constraints} 22:1--29.

\bibitem[{Joshi, Laurent, \protect\BIBand{} Bresson(2019)}]{joshi2019efficient}
Joshi CK, Laurent T, Bresson X, 2019 \emph{An efficient graph convolutional network technique for the travelling salesman problem} \urlprefix\url{http://arxiv.org/abs/1906.01227}.

\bibitem[{Kalcsics \protect\BIBand{} R{\'i}os-Mercado(2019)}]{kalcsics2019}
Kalcsics J, R{\'i}os-Mercado RZ, 2019 \emph{Districting Problems}, 705--743 (Springer International Publishing).

\bibitem[{Kingma \protect\BIBand{} Ba(2014)}]{kingma2014adam}
Kingma DP, Ba J, 2014 \emph{Adam: A method for stochastic optimization} \urlprefix\url{http://arxiv.org/abs/1412.6980}.

\bibitem[{Kipf \protect\BIBand{} Welling(2017)}]{kipf2016semi}
Kipf TN, Welling M, 2017 \emph{Semi-supervised classification with graph convolutional networks}. \emph{International Conference on Learning Representation}.

\bibitem[{Kool, van Hoof, \protect\BIBand{} Welling(2018)}]{kool2018attention}
Kool W, van Hoof H, Welling M, 2018 \emph{Attention, learn to solve routing problems!} \emph{International Conference on Learning Representations}.

\bibitem[{Kou, Golden, \protect\BIBand{} Bertazzi(2024)}]{Kou2024}
Kou S, Golden B, Bertazzi L, 2024 \emph{{An improved model for estimating optimal VRP solution values}}. \emph{Optimization Letters} 18(3):697--703.

\bibitem[{Kou, Golden, \protect\BIBand{} Poikonen(2022)}]{Kou2022a}
Kou S, Golden B, Poikonen S, 2022 \emph{{Optimal TSP tour length estimation using standard deviation as a predictor}}. \emph{Computers \& Operations Research} 148:105993.

\bibitem[{Kovacs et~al.(2014)Kovacs, Golden, Hartl, \protect\BIBand{} Parragh}]{Kovacs2014}
Kovacs A, Golden B, Hartl R, Parragh S, 2014 \emph{{Vehicle routing problems in which consistency considerations are important: A survey}}. \emph{Networks} 64(3):192--213.

\bibitem[{Kwon, Golden, \protect\BIBand{} Wasil(1995)}]{Kwon1995}
Kwon O, Golden B, Wasil E, 1995 \emph{{Estimating the length of the optimal TSP tour: An empirical study using regression and neural networks}}. \emph{Computers {\&} Operations Research} 22(10):1039--1046.

\bibitem[{Kwon et~al.(2020)Kwon, Choo, Kim, Yoon, Gwon, \protect\BIBand{} Min}]{kwon2020pomo}
Kwon YD, Choo J, Kim B, Yoon I, Gwon Y, Min S, 2020 \emph{{POMO: Policy optimization with multiple optima for reinforcement learning}}. \emph{Proceedings of the 34th International Conference on Neural Information Processing Systems}.

\bibitem[{Lei, Laporte, \protect\BIBand{} Guo(2012)}]{Lei2012}
Lei H, Laporte G, Guo B, 2012 \emph{{Districting for routing with stochastic customers}}. \emph{EURO Journal on Transportation and Logistics} 1(1):67--85.

\bibitem[{Lei et~al.(2015)Lei, Laporte, Liu, \protect\BIBand{} Zhang}]{Lei2015}
Lei H, Laporte G, Liu Y, Zhang T, 2015 \emph{{Dynamic design of sales territories}}. \emph{Computers {\&} Operations Research} 56:84--92.

\bibitem[{Levy-Kramer(2018)}]{Levy-Kramer_k-means-constrained_2018}
Levy-Kramer J, 2018 \emph{{k-means-constrained}}. \urlprefix\url{https://github.com/joshlk/k-means-constrained}.

\bibitem[{Li et~al.(2016)Li, Tarlow, Brockschmidt, \protect\BIBand{} Zemel}]{li2015gated}
Li Y, Tarlow D, Brockschmidt M, Zemel R, 2016 \emph{Gated graph sequence neural networks}. \emph{International Conference on Learning Representation}.

\bibitem[{Lin \protect\BIBand{} Kernighan(1973)}]{LinKernighan1973}
Lin S, Kernighan BW, 1973 \emph{An effective heuristic algorithm for the traveling-salesman problem}. \emph{Operations Research} 21(2):498–516.

\bibitem[{Louren{\c{c}}o, Martin, \protect\BIBand{} St{\"u}tzle(2019)}]{Lourenco2010ILS}
Louren{\c{c}}o HR, Martin OC, St{\"u}tzle T, 2019 \emph{Iterated local search: Framework and applications}. \emph{Handbook of Metaheuristics} 129--168.

\bibitem[{Miyazawa et~al.(2020)Miyazawa, Moura, Ota, \protect\BIBand{} Wakabayashi}]{miyazawa2020flow}
Miyazawa FK, Moura PFS, Ota MJ, Wakabayashi Y, 2020 \emph{Cut and flow formulations for the balanced connected k-partition problem}. Ba{\"i}ou M, Gendron B, G{\"u}nl{\"u}k O, Mahjoub AR, eds., \emph{Combinatorial Optimization}, 128--139 (Springer International Publishing).

\bibitem[{Novaes, de~Cursi, \protect\BIBand{} Graciolli(2000)}]{Novaes2000}
Novaes A, de~Cursi J, Graciolli O, 2000 \emph{{A continuous approach to the design of physical distribution systems}}. \emph{Computers {\&} Operations Research} 27(9):877--893.

\bibitem[{Park(2018)}]{CensusData2018}
Park N, 2018 \emph{Middle super output area population estimates - mid-2018: Sape21dt3a edition}. \url{https://www.ons.gov.uk/peoplepopulationandcommunity/populationandmigration/populationestimates/datasets/middlesuperoutputareamidyearpopulationestimates}, accessed on November 27, 2021.

\bibitem[{Scarselli et~al.(2009)Scarselli, Gori, Tsoi, Hagenbuchner, \protect\BIBand{} Monfardini}]{Scarselli2009}
Scarselli F, Gori M, Tsoi AC, Hagenbuchner M, Monfardini G, 2009 \emph{The graph neural network model}. \emph{IEEE Transactions on Neural Networks} 20(1):61--80.

\bibitem[{Sun \protect\BIBand{} Yang(2023)}]{sun2023difusco}
Sun Z, Yang Y, 2023 \emph{Difusco: Graph-based diffusion solvers for combinatorial optimization}. \emph{Advances in neural information processing systems} 36:3706--3731.

\bibitem[{Uchoa et~al.(2017)Uchoa, Pecin, Pessoa, Poggi, Vidal, \protect\BIBand{} Subramanian}]{Uchoa2017}
Uchoa E, Pecin D, Pessoa A, Poggi M, Vidal T, Subramanian A, 2017 \emph{New benchmark instances for the capacitated vehicle routing problem}. \emph{European Journal of Operational Research} 257(3):845--858.

\bibitem[{Varol, {\"{O}}zener, \protect\BIBand{} Albey(2024)}]{Varol2024}
Varol T, {\"{O}}zener O, Albey E, 2024 \emph{{Neural network estimators for optimal tour lengths of traveling salesperson problem instances with arbitrary node distributions}}. \emph{Transportation Science} 58(1):45--66.

\bibitem[{Veli\v{c}kovi\'{c} et~al.(2018)Veli\v{c}kovi\'{c}, Cucurull, Casanova, Romero, Li{\`{o}}, \protect\BIBand{} Bengio}]{velickovic2018graph}
Veli\v{c}kovi\'{c} P, Cucurull G, Casanova A, Romero A, Li{\`{o}} P, Bengio Y, 2018 \emph{Graph attention networks}. \emph{International Conference on Learning Representations}.

\bibitem[{Verweij et~al.(2003)Verweij, Ahmed, Kleywegt, Nemhauser, \protect\BIBand{} Shapiro}]{Verweij2003}
Verweij B, Ahmed S, Kleywegt AJ, Nemhauser G, Shapiro A, 2003 \emph{The sample average approximation method applied to stochastic routing problems: A computational study}. \emph{Computational Optimization and Applications} 24(2–3):289–333.

\bibitem[{Wang(2019)}]{wang2019deep}
Wang MY, 2019 \emph{Deep graph library: Towards efficient and scalable deep learning on graphs}. \emph{ICLR workshop on representation learning on graphs and manifolds}.

\bibitem[{Webster(2013)}]{Webster2013}
Webster GR, 2013 \emph{{Reflections on current criteria to evaluate redistricting plans}}. \emph{Political Geography} 32(1):3--14.

\bibitem[{Young(1988)}]{Young_1988}
Young HP, 1988 \emph{Measuring the compactness of legislative districts}. \emph{Legislative Studies Quarterly} 13(1):105.

\bibitem[{Zhong, Hall, \protect\BIBand{} Dessouky(2007)}]{Zhong2007}
Zhong H, Hall RW, Dessouky M, 2007 \emph{{Territory planning and vehicle dispatching with driver learning}}. \emph{Transportation Science} 41(1):74--89.

\bibitem[{Zoltners \protect\BIBand{} Sinha(2005)}]{Zoltners2005}
Zoltners AA, Sinha P, 2005 \emph{{Sales territory design: Thirty years of modeling and implementation}}. \emph{Marketing Science} 24(3):313--331.

\end{thebibliography}
\end{document}